%% file: SMO_Paper2023.tex
\documentclass[twocolumn]{svjour3}          

\usepackage[colorlinks,citecolor=black,linkcolor=black,urlcolor=black]{hyperref}
\usepackage{amsmath}
\usepackage{amssymb}
\usepackage{bm}
\usepackage{graphicx}
\usepackage{booktabs}
\usepackage{appendix}
\usepackage{doi}
\usepackage{url}
\usepackage{hyphenat}
\usepackage{float}
\usepackage[natbib=true,citestyle=authoryear,backend=bibtex,url=false,maxcitenames=2]{biblatex}
\usepackage{stmaryrd} 
\usepackage{mathtools} 
\usepackage{subfigure}
\usepackage{comment}
\usepackage[percent]{overpic}
\usepackage{multirow} 
\usepackage{booktabs} 
\usepackage{makecell} 

\usepackage{microtype}
\usepackage{float}
\usepackage{amsopn}

\DeclareMathOperator{\atan2}{atan2}

\let\cite\textcite
\journalname{Structural and Multidisciplinary Optimization}

\begin{document}

\title{Concurrent level set and fiber orientation optimization of composite structures}

\author{M.Mokhtarzadeh$^{1}$, F López Jiménez$^{2}$ and K. Maute$^{3}$}

\institute{
$^1${Aerospace Mechanics Research Center, Department of Aerospace Engineering Sciences},\\
{University of Colorado Boulder}, {3775 Discovery Dr}, {Boulder}, {80309-0429}, {CO}, {USA}\\
\email{momo1035@colorado.edu}\\
$^2${Aerospace Mechanics Research Center, Department of Aerospace Engineering Sciences},\\
{University of Colorado Boulder}, {3775 Discovery Dr}, {Boulder}, {80309-0429}, {CO}, {USA}\\
\email{Francisco.LopezJimenez@colorado.edu}\\
$^3${Aerospace Mechanics Research Center, Department of Aerospace Engineering Sciences},\\
{University of Colorado Boulder}, {3775 Discovery Dr}, {Boulder}, {80309-0429}, {CO}, {USA}\\
\email{k.maute@colorado.edu} \\
}

\date{Received: date / Accepted: date}

\maketitle

\abstract{By adjusting both the structural shape and fiber orientation, this research aims to optimize the design of Fiber Reinforced Composite (FRC) structures. 
The structural geometry is represented by a level set function, which is approximated by quadratic B-spline functions. 
The fiber orientation field is parameterized with quadratic/cubic B-splines on hierarchically refined meshes. Different levels 
for B-spline mesh refinement for the level set and fiber orientation fields are studied to resolve geometric features and to obtain a smooth fiber layout. 
To facilitate FRC manufacturing, the parallel alignment, and smoothness of fiber paths are enforced by introducing penalty terms referred to as "misalignment penalty and 
curvature penalty", which are incorporated into the optimization process. A geometric interpretation of the penalties is provided. 
The material behavior of the FRCs is modeled by the Mori-Tanaka homogenization scheme and the macroscopic structure response is 
modeled by linear elasticity under static mutiloading conditions. The Governing equations are discretized by a Heaviside-enriched eXtended IsoGeometric Analysis (XIGA) to avoid the need to 
generate conformal meshes. Instabilities in XIGA are mitigated by the facet-oriented ghost stabilization technique. This work considers mass and strain energy in the 
formulation of the optimization objective, along with misalignment and curvature penalties and additional regularization terms. Constraints are imposed on the volume of the structure. 
The resulting optimization problems are solved by a gradient-based algorithm. The design sensitivities are computed by the adjoint method. Numerical examples demonstrate with 
two-dimensional and three-dimensional configurations that the proposed method is efficient in simultaneously optimizing the macroscopic shape and the fiber layout while improving 
manufacturability by promoting parallel and smooth fiber paths. }

\maketitle
\input{symbolss.tex}
\input{Intro.tex}

\input{Analysis.tex}
\input{Opt.tex}

\input{Examples.tex}
















\section{Conclusions}\label{sectionConclusion}

In this study, we present a concurrent optimization approach for the topology and fiber orientation of 
Fiber Reinforced Composite (FRC) structures. 
The structural shape is represented by a level set function approximated through quadratic B-splines, 
while the fiber orientation fields are discretized with linear, quadratic, and cubic B-splines with coarser meshes. 
Employing level sets for geometry definition ensures precise boundary definition 
while utilizing higher-order coarse B-splines for the discretization of the fiber orientation and level set promotes smoother designs. 
The analysis employs an eXtended Finite Element Method (XFEM) framework, adopting a generalized Heaviside enrichment and facet-oriented ghost stabilization for 
improved stability and robustness. 
FRCs' material behavior is modeled with linear elasticity, with the elasticity tensor derived from Mori-Tanaka homogenization.

Novel penalty terms are introduced to promote parallel and smooth fiber paths. 
These terms anisotropically constrain the fiber orientation field's first-order derivative, 
and are integrated into the optimization problem formulation via the penalty terms. 
The geometric interpretation of these penalties is explored by defining optimization problems to assess their impact on a hypothetical fiber 
path and their efficacy in creating parallel and smooth paths.

This study considered compliance minimization problems and examined 
the interplay between geometry and fiber orientation 
for different penalty combinations. The influence of the discretization of the fiber orientation field in terms of polynomial order and
B-spline mesh size and its influence on the optimized structure were studied. 

Numerical experiments with various penalty combinations demonstrated the 
effectiveness of the parallel misalignment penalty in aligning fiber paths in both 
2D and 3D settings. 
Simultaneously, the curvature penalties were found to be efficient in controlling bending and waviness 
in the fiber paths. 
Incorporating the misalignment and curvature penalties in the optimization problem formulation led to an 
increase in optimized strain energy values. 
The strain energy increase was more pronounced with the addition of more penalty terms, 
indicating that increased parallelism and smoothness 
in the fiber paths lead to higher strain energy values. 

The numerical examples suggested that the polynomial order of fiber orientation did not 
significantly influence strain energy, yielding similar outcomes across all polynomial orders. 
While linear B-splines failed to produce smooth curvature variations, this lack of smoothness did 
not impact the smoothness of the streamline visualization of the fiber path. 
Coarser discretization of the fiber orientation fields leads to a smoother fiber layout by reducing localized features. However, this smoothness 
comes at the cost of increased strain energy values.

Optimized fiber paths, characterized by their parallel alignment and smoothness, simplify the 
post-processing efforts and narrow the difference between the optimized layout and the post-processed, 
uniformly spaced fiber path.

While the numerical experiments render the proposed method promising, a fundamental problem remains to be addressed. 
Although $\pi$-periodicity issues did not arise in the numerical examples studied in this paper, a similar issue was encountered in the shear-dominated examples, see Section \ref{example1} and \ref{example4}. 
In a state of pure shear, the fiber angles $\frac{\pi}{4}$ and $\frac{3\pi}{4}$ result in the same stiffness. The optimization process may create regions where the fiber orientation is $\frac{\pi}{4}$ and nearby regions where the fiber orientation is $\frac{3\pi}{4}$.  The continuous interpolation of the fiber orientation field used in the proposed method necessitates transition regions where the fiber angles vary between  $\frac{\pi}{4}$ and $\frac{3\pi}{4}$. Thus, in these transition regions, the angle is not optimal, affecting the performance of the design. A similar issue may arise for general loading scenarios due to the $\pi$-periodicity.  In future studies the proposed method needs to be enhanced to avoid the need for these transition regions by, for example, allowing for a locally discontinuous fiber orientation field.


\section*{Acknowledgments}
The authors acknowledge the support from the Air Force Office of Scientific Research (AFOSR) under Grant FA9550-20-1-0306, with Dr. Byung-Lip Lee as the program manager. 
The views expressed are those of the authors and do not reflect the official policy or position of the AFOSR or the U.S. Government.

The authors would like to thank Dr. Narasimha Boddeti for providing the code for generating continuous fiber paths.

\section*{Compliance with ethical standards}
The opinions and conclusions presented in this paper are those of the authors and do not necessarily reflect the views of the sponsoring organizations.

\section*{Conflicts of Interest} 
On behalf of all authors, the corresponding author states that there is no conflict of interest. 

\section*{Replication of Results} 
Upon request, the authors will provide the full set of input parameters and meshes for each validation problem and topology optimization problem presented in the paper.


\begin{appendices}
\end{appendices}


\printbibliography
\end{document}

%% file: symbolss.tex
\newcommand{\strainEnergy}{\mathcal{F}}

%% file: Intro.tex
\section{Introduction}\label{sectionIntroduction}
Fiber Reinforced Composites (FRC) are materials that are
made up of a matrix (the continuous phase) and fibers (the dispersed phase).
The fibers are embedded in the matrix and provide the composite material
with enhanced mechanical properties, such as increased stiffness, strength, 
and fatigue resistance. Therefore, these materials exhibit a superior stiffness-to-weight ratio compared to conventional isotropic homogeneous materials.  

In the early stages of FRC design, 
the primary emphasis was on determining the orientation of fibers for given structural shapes like beams or plates, 
as noted in \cite{nikbakt2018review}. These designs typically assumed a constant fiber orientation throughout each ply.
However, advancements in composite manufacturing techniques, such as continuous fiber fused filament fabrication (CF4) \citep{wang2021load},
allow for spatially varying fiber angles, increasing structural performance if the local angles are chosen appropriately. 
Moreover, if the structural shape can be designed alongside the fiber orientation further performance boost can be expected. 
This paper introduces a design optimization framework for the simultaneous optimization of 
the structural geometry and fiber orientation in FRC structures.

To optimize the structural shape, this paper uses Topology Optimization (TO), which is a reliable method for finding the optimal layout of a structure. 
It provides a systematic way to alter the topology and shape of a structure and effectively 
eliminate unnecessary material. TO does not require a close-to-optimal design to initialize the optimization process and can significantly reduce its weight and cost.

In this work, design optimization is used to determine not only the topology of the structure but also the fiber orientation to achieve a desired set of mechanical properties.
There have been several attempts to optimize fiber orientation and topology simultaneously in the literature, 
most of which are summarized in the review paper by \cite{2022gandhireview}. 
Concurrent fiber and topology optimization involves two main components: 
representation of structural geometry and parameterization of fibers. The past work differs mainly in the latter component.

To represent the geometry of the structure, most of the methods proposed in the literature 
employ the density method which represents the material distribution as a volume fraction field, ranging from 0 (indicating void) to 1 (indicating solid).
The material properties are typically interpolated using the Solid Isotropic Material with Penalization (SIMP) approach, 
often in conjunction with filtering and projection techniques. These additional steps are crucial for controlling feature 
size and accurately delineating the solid-void boundary interface, as exemplified in \cite{sigmund2013topology}. Conversely, the Level Set Method (LSM) 
addresses the challenge of defining this interface, offering a precise description of structural geometry throughout the optimization process and eliminating 
the need for a material interpolation scheme, as discussed by \cite{van2013level}. In this work, the LSM is selected to describe the geometry in the current study.

Methods for concurrent optimization of FRCs vary by the type of design parameters which may include fiber orientation, fiber volume fraction, 
and other fiber properties. The focus of this paper is on FRC structures with designable spatially varying fiber orientation. 
Depending on the manufacturing techniques, different levels of spatial variability of the fiber orientation can be realized. For example, CF4 allows for a continuously varying fiber orientation, 
whereas other manufacturing techniques such as Automated Tape Layout (ATL) only allow for a discrete set of fiber orientation values.
In the context of continuous fiber orientation optimization (CFOA) for FRCs,
 maintaining a continuous and smooth fiber trajectory is crucial. 
 This continuity is imperative for two primary reasons: 
 firstly, interruptions in fiber alignment can lead to stress concentrations, 
 weakening the composite structure. Secondly, 
 a continuous, smooth, and parallel fiber path is a prerequisite for most manufacturing techniques.

In discrete parameterization approaches to fiber orientation optimization, 
the orientation is treated as a discrete variable, limited to a predefined set of angles. 
This approach, known as Discrete Material Optimization (DMO), was initially introduced in the work of
\cite{Stegmann2005}. Subsequent advancements, akin to those in \cite{Kiyono2017}, 
expanded the methodology by augmenting the discrete angle set with angles selected from a normal distribution.
 However, DMO's inherent constraint of limiting fiber 
orientations to a predefined set makes it less suitable for manufacturing techniques like CF4. 
Such techniques offer broader design freedoms, which DMO does not fully exploit due to its discrete nature.

To fully leverage the design capabilities of Additive Manufacturing (AM) methods such as the CF4 method, 
it is essential to allow for continuous variation in fiber orientation. 
The initial approach in this field, the stress/strain-based design, was introduced by 
\cite{Pedersen1989} and further developed by \cite{Gea}. This method aligns fiber orientations with the principal 
stress or strain directions. This method is not applicable to problems with multiple load cases and does 
not formally consider fiber orientation as an optimization variable which limits its broader application in an optimization
algorithm.

Treating fiber angle as a continuous design variable in the context of 
laminate composites was initially introduced by \cite{Bruyneel2002}, allowing for a range of fiber 
angles from $[0, \pi]$. In manufacturing, the application of spatially continuous design spaces often necessitates a 
smooth fiber layout. To address this, several studies, such as \cite{Papapetrou2020}, \cite{Brampton2015}, 
\cite{bruyneel2013modified}, and \cite{fernandez2019optimal}, have aligned fibers along iso-contours of a 
level set field or streamlines. This method produces parallel, equidistant fibers, 
simplifying the post-processing steps required for their manufacturing. 
However, this method restricts the design space, potentially leading to sub-optimal designs, 
as discussed in \cite{tian2021parametric}.

Feature-based parametrization is another technique for obtaining spatially continuous fiber orientation.
This approach involves designing bars reinforced with continuous fibers aligned along the bars, see \cite{smith2021topology}, 
and \cite{greifenstein2023efficient}. 
While effective, this method also imposes limitations on the design space and links fiber orientation 
with the structure's geometry and does not take full advantage of the design freedom offered by AM.

The isoparametric transformation method is a continuous fiber orientation optimization approach 
that extends the principles of the density method. Introduced by \cite{Nomura2015a}, 
it converts the fiber angle into a Cartesian vector with two independent components. These components 
along with the density variable, are used in the SIMP method for material property interpolation 
of anisotropic materials. This method, although effective in creating a smooth fiber field through 
filtering and projection, doubles the design variables and requires a constraint for maintaining the 
unit norm of the angle vector. It also faces challenges similar to the density method, like defining interfaces and boundaries.
Follow-up studies by \cite{Kim2020}, \cite{jung2022inverse} and \cite{smith2022topology} have further developed \cite{Nomura2015a}'s method. 
To ensure fiber path smoothness, these methods typically adjust filter radius and constrain design variable derivatives, 
as explored by \cite{greifenstein2016simultaneous}. 
Additionally, some researchers, including \cite{Papapetrou2020}, \cite{boddeti2020optimal}, \cite{fedulov2021optimization}, 
and \cite{fernandes2021experimental}, have investigated post-processing techniques to achieve smooth and parallel fiber 
layouts from optimized orientation fields to bridge the gap between the optimized fiber orientation field and the final manufacturable structure's fiber layout.

A few studies have integrated constraints within the optimization process to ensure smooth and parallel fiber layouts. \cite{tian2019optimization} 
employed overlap and curvature constraints in a fiber orientation field defined through Shepard interpolation, specifying areas for 
parallel and curvatures-constrained fiber orientation. However, 
this approach is limited by its reliance on a predetermined geometry and a 2D, 
mesh-based formulation for fiber orientation constraints, which complicates its extension to 
evolving topologies and 3D problems. Similarly, \cite{fernandez2019optimal} 
introduced a method that integrates manufacturing constraints into the optimization process, 
but it is also based on a fixed geometry assumption.

A design optimization methodology is presented in this paper to address the shortcomings discussed above. 
The structure's geometry is represented using level sets, ensuring a clear interface definition throughout the optimization.
The fiber orientation is treated as a continuously variable field for greater design flexibility. 
In order to achieve this continuous interpretation, higher-order B-splines with adjustable 
refinement levels are employed to parametrize the fiber orientation field. B-spline functions promote smooth designs, prevent the 
appearance of spurious features, and eliminate the need for further filtering techniques, as discussed by Noël et al. (2020).
The fiber orientation field is filtered implicitly by adjusting the refinement level and polynomial order of the B-spline.

Misalignment and curvature penalties are introduced to promote parallel and smooth fiber alignment. 
These penalties control the first-order derivative of the fiber orientation field in an anisotropic manner and
are incorporated into the formulation of the objective function.
A detailed geometric interpretation is provided to explain how these penalty terms facilitate the creation of smooth and parallel fiber layouts. 
While some post-processing of the optimized structure is necessary, 
these penalties reduce the effort required in this phase by generating a smoother and more parallel fiber layout,  
effectively reducing the gap between the optimized fiber orientation field and the final, manufactured structure.

In this study, the structural behavior is modeled using static linear elasticity, 
while the microscopic behavior of the FRC is captured through the Mori-Tanaka homogenization scheme \citep{mori1973average}.
The weak form of the governing equations is discretized using eXtended IsoGeometric Analysis (XIGA). 
Quadratic and cubic B-splines are employed for the level set, fiber orientation field, and state variables. 
Enrichment in XIGA is achieved using a generalized Heaviside approach, 
and the Nitsche method \citep{Nitsche1971} is employed to weakly enforce Dirichlet boundary conditions. Furthermore, 
to tackle numerical instabilities in small subdomains, which arise from the immersed nature of XIGA, 
facet-oriented ghost stabilization is employed.
The resulting optimization problems are solved by the Globally Convergent Method of Moving Asymptotes (GCMMA). 
The gradients of the objective and constraint functions are computed by the adjoint method.

The remainder of this paper is organized as follows:
Section \ref{macromicro} discusses the level set formulation for geometry representation, 
the B-spline parameterization of the fiber orientation field. 
Section \ref{analysis} provides a brief description of XIGA and its building blocks as well as the weak form of the governing equations. 
An overview of the formulation of the optimization problem and the misalignment and curvature penalties is provided in 
Section \ref{optimzation}. Section \ref{sectionExamples} presents numerical examples where the effect of misalignment and smoothing penalties 
and parameterization of the fiber orientation field on the optimized design is discussed. 
Section \ref{sectionConclusion} summarizes the findings of this study and draws conclusions about the developed method.

%% file: Analysis.tex
\section{Design variable representation}\label{macromicro}
In this study, design variables for concurrent optimization of the FRC structures include geometry
and fiber orientation. Section 2.1 details the parametrization of geometry, and Section 2.2 discusses 
the parametrization of fiber orientation.
 
\subsection{Geometry representation}\label{geomtryrepresnt}
The shape of the structure and material interfaces are described by 
a level set function (LSF), \( \phi \); see for example \cite{van2013level} and references therein. The interface and 
external boundaries are defined as the zero iso-level of the LSF.
The LSF is a scalar function that discriminates between the two material domains, \( \Omega_1 \) and \( \Omega_2 \), by assigning positive 
and negative values respectively, with \( \Gamma_{12} \) representing the boundary. 
Note that one of these domains may represent void. Formally, the LSF at a spatial point with 
coordinate \( \mathbf{x} \), 
within the computational domain is given by:

\begin{equation}
	\phi(\mathbf{x}) = 
	\begin{cases}
	< 0, & \text{for all } \mathbf{x} \in \Omega_1, \\
	> 0, & \text{for all } \mathbf{x} \in \Omega_2, \\
	= 0, & \text{for all } \mathbf{x} \in \Gamma_{12}.
	\end{cases}
	\label{eq:LSF_definition}
\end{equation}

\begin{figure}[!h]
	\centering
	\includegraphics[scale=1.2]{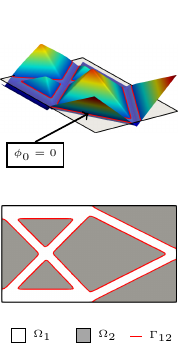}
	\caption{Geometric description of solid (\( \Omega_1 \))/void (\( \Omega_2 \)) design domain using the LSF.}
	\label{fig_LS_geometry}
\end{figure}

Figure \ref{fig_LS_geometry} illustrates the LSF's application in defining the geometry of a two-phase solid/void design domain, with \( \Omega_1 \) indicating the 
solid structure, \( \Omega_2 \) the void, and \( \Gamma_{12} \) the interface between them.

The LSF is discretized on a computational mesh which may be independent of the state variable mesh. In this work, the LSF is approximated using B-spline basis functions 
\( B_k(\mathbf{x}) \) and their corresponding coefficients \( \phi_k \) as:

\begin{equation}
\phi(\mathbf{x}) \approx \phi^h(\mathbf{x}) = \sum_{k} B_k(\mathbf{x}) \phi_k,
\label{eq:LSF_discretization}
\end{equation}

\noindent where \( \phi^h(\mathbf{x}) \) denotes the discretized LSF. 

In the work of \cite{sethian2000structural}, optimization with the level set method involves evolving the Level Set Function (LSF) 
through the Hamilton-Jacobi equation. In contrast, the current study follows the work of \cite{Noel2022} and defines the LSF coefficients explicitly as functions of the design variables, 
represented by vector $\bm{s}$ of length $N_s$. These coefficients are then updated via a gradient-based algorithm employing shape sensitivities.
This approach simplifies managing multiple constraints.

\subsection{Fiber orientation}\label{FO}

This work considers FRCs in 
two dimensions with one spatially varying fiber direction and three dimensions with two fiber angles.
 Figure \ref{fig_angles} illustrates the geometrical definitions of the angles for fiber orientation
 (\( \theta_{xy} \) , \( \theta_z \))  
. The tangent vector \( \mathbf{t} \) along a fiber's path is 
projected onto the x-y plane, establishing \( \theta_{xy} \) as the first angle. 
Conversely, \( \theta_z \) is defined as the angle between the fiber's tangent and its 
projection onto the x-y plane. Utilizing these definitions, the fiber's tangent vector is expressed in two and three dimensions as follows:

\begin{figure}[!h]\centering
	\includegraphics[scale=1.0]{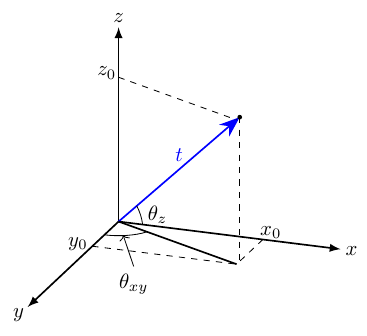}
	\caption{Fiber angle definitions in 3D ($\theta_{xy}$,$\theta_z$).}
	\label{fig_angles}
\end{figure}

\begin{equation}
\mathbf{t}_{2D} = \begin{bmatrix}
\cos(\theta_{xy}) \\
\sin(\theta_{xy}) 
\end{bmatrix}, \quad
\mathbf{t}_{3D} = \begin{bmatrix}
\cos(\theta_{xy})\cos(\theta_z) \\
\sin(\theta_{xy})\cos(\theta_z) \\
\sin(\theta_z)
\end{bmatrix}.
\label{eqTangentVector}
\end{equation}

In this study, we focus on continuously varying fiber orientation and to fully leverage the design potential offered by advanced manufacturing processes, 
we adopt a field-based approach. This approach represents fiber orientation as a continuous field in the computational domain, 
unlike the iso-contours approach of \cite{bruyneel2013modified} or feature-mapping approach of \cite{greifenstein2023efficient} where fibers are aligned with bars. 
These field-based approaches can be viewed as an extension of the density method. 

Among field-based methods, \cite{Nomura2015a} introduced a new variation, converting the fiber orientation field into a Cartesian vector field, addressing the 
$2\pi$ periodicity issue where angles like $0$ and $2\pi$ represent the same orientation.
In our approach, where the fiber angle is directly utilized as a design variable, we extend the bounds beyond the conventional $[0,2\pi]$ range.
Notably, as discussed in our numerical examples in Section \ref{sectionExamples}, the fiber orientation does not reach the prescribed upper and lower bounds. 
It's important to note that the Cartesian transformation doesn't resolve the $\pi$-periodicity in the elastic properties of materials, resulting in identical material properties for angles 
differing by $\pi$. Employing a continuous interpolation method may create a transition zone between angles associated with $\pi$-periodicity,
resulting in sub-optimal designs.

Most field-based methods, including those in \cite{almeida2023concurrent} and \cite{smith2022topology}, 
are typically paired with the density method and adopt the same interpolation for fiber orientation as the density variables. 
Other parametrization approaches exist, as shown in \cite{tian2019optimization}, where fiber orientation is constant within each element in the Finite Element Method (FEM) mesh.
In this work, we utilize B-spline interpolation for the fiber orientation field, 
chosen for their smoothness. 

For 2D problems, only one B-spline discretization suffices to describe the tangent vector
field; see Equation (\ref{eqTangentVector}). 3D 
problems necessitate two independent B-spline discretizations—one for the $\theta_{xy}$ and another for the $\theta_z$. 
The fiber orientation fields are discretized on a mesh that may differ in polynomial order and refinement level from those used for 
state variables and geometry. This discretization of the fiber orientation field, $\theta$, is achieved through B-spline basis functions as:

\begin{equation}
\theta(\mathbf{x}) \approx \theta^h(\mathbf{x}) = \sum_{k} B_k(\mathbf{x}) \theta_k,
\label{eqFiberBspline1}
\end{equation}

\noindent where \( B_k(\mathbf{x}) \) denotes the B-spline basis functions and \( \theta_k \) the coefficients of the approximating function.
In Equation (\ref{eqFiberBspline1}), $\theta$ represents either $\theta_{xy}$ or $\theta_z$. 
Analogous to the level set functions, these B-spline coefficients are updated 
via a gradient-based algorithm
using material parameter sensitivities.

\begin{figure}[t]
	\centering
	\includegraphics[scale=1]{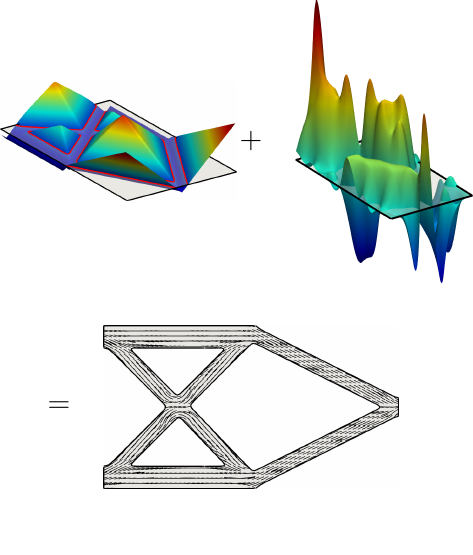}
	\caption{Illustration of level set function for a truss structure (top left), the B-spline surface parameterizing the fiber orientation(top right), 
	and the resulting fiber orientation and geometry (bottom).}
	\label{fig_bspline_angles}
\end{figure}

Figure \ref{fig_bspline_angles} illustrates the overlay of the geometry with the fiber orientation. 
The top left and top right sections of this figure depict B-spline surfaces, representing the level set function and fiber orientation, respectively. 
The lower part of the figure illustrates the superposition of these two B-spline surfaces, 
showcasing the resulting structure's geometry and its fiber orientation.

\section{Structural Analysis}\label{analysis}
In the present work, the eXtended IsoGeometric Analysis (XIGA) methodology is employed, 
which builds upon the traditional eXtended Finite Element Method (XFEM) framework by incorporating B-spline functions as the basis, 
as discussed in \cite{Noel2022}. XFEM itself is a variant of the classical Finite Element Method (FEM), 
specifically tailored to immersed geometry descriptions. 
XFEM can be seamlessly combined with level set topology optimization
to manage evolving design interfaces, eliminating the need for generating conformal meshes. This integration maintains 
the crisp definition of the interface, as represented in the level set, within the physical FEM model.
This extension to the FEM involves the augmentation of additional basis functions and Degrees Of Freedom (DOFs) — an approach termed `enrichment` — 
to capture the physical response near interfaces and boundaries. The following subsections will summarize the enrichment strategy 
 and present the weak formulation of the governing equations.

\subsection{Enrichment strategy}
This work utilizes the generalized Heaviside enrichment strategy described in \cite{2020NoelEtAl}. 
This enrichment approach accommodates a variety of material phases, intersection configurations, 
and basis function supports.

\begin{figure}[t]\centering
	\includegraphics[scale=0.8]{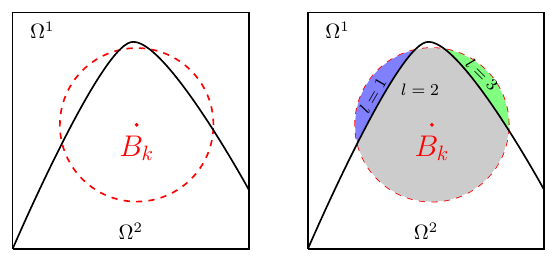}
	\caption{Enrichment strategy in XIGA.}
	\label{fig_enrichment}
\end{figure}

Consider the configuration in Fig.~\ref{fig_enrichment} which shows a region covered by two material phases ($\Omega_1$ and $\Omega_2$) 
and the support of the basis function, $B_k$,  indicated by dashed red lines. 
The basis support consists of three distinct subregions. To accurately represent the physical response in these subregions without spurious coupling,
the same basis function is weighted by three different independent coefficients, i.e., Degrees of Freedom (DOFs).
Phase 1 occupies subregions $l=1$ and $l=3$, which are topologically disconnected, while phase 2 occupies the subregion $l=2$. The set of all topologically disconnected 
subregions are denoted by $\{ \Omega_k^{\ell} \}_{\ell=1}^{L_k}$, where $L_k$ is the total number of these subregions. 
In general, the $i^{th}$ component of the discretized vector-valued state variable field, denoted as $\mathbf{u}_i$, which corresponds to displacement in this study, is expressed as:

\begin{equation}
    \mathbf{u}_i(\mathbf{x}) = \sum_{k=1}^{K} \sum_{\ell=1}^{L_k} \varphi_k^{\ell}(\mathbf{x})\, B_k(\mathbf{x})\, u_{i,k}^{\ell},
\end{equation}

\noindent where the total number of background basis functions is denoted by $K$.
The coefficients ${u}^{l}_{i,k}$ are the DOFs associated with the basis function $B_k$ and subregions $\Omega_{k}^{\ell}$.

The indicator function $\varphi_k^{\ell}(\mathbf{x})$ is a binary-valued function which indicates membership of $\mathbf{x}$  in $\Omega_k^{\ell}$ and is defined as:

\begin{equation}
\varphi_k^{\ell}(\mathbf{x}) = I_{\Omega_k^{\ell}}(\mathbf{x}) = \left\{
\begin{array}{ll}
1, & \displaystyle \mbox{if}\ \mathbf{x} \in \Omega_k^{\ell},\\[5pt]
0, & \displaystyle \mbox{otherwise}.
\end{array}
\right.
\label{Eq_enrichment}
\end{equation}

\subsection{Governing equations}\label{governeqn}

The enrichment strategy is applied to the weak form of the governing equations, 
which models the static linear elastic response of a structure. 
The residual equations are broken down into four separate components:
\begin{equation}\label{Eqweakform}
\mathcal R = \mathcal R ^{ U }+ \mathcal R ^{ D }+ \mathcal R ^{ G } + \mathcal R ^{ S } =0. 
\end{equation}

The bulk contribution and Neumann boundary contribution are included in $ \mathcal R ^{ U }$, the Dirichlet boundary condition represented by $\mathcal R ^{ D }$, 
the facet-oriented ghost stabilization term added through $ \mathcal R ^{ G }$, and $\mathcal R ^{ S }$ 
represents a stabilization term to suppress rigid body motions. 

Assuming linear elasticity without body forces under static loading, the bulk and Neuman contributions are as follows:

\begin{equation}
\mathcal R ^{ U }= \int_{\Omega^m}  \bm{\varepsilon}(\delta \mathbf  u) : \bm  \sigma d \Omega-\int_{\Gamma_{N}} \delta \mathbf  u \cdot \mathbf  f _N d \Gamma,
\end{equation}

\noindent where $\mathbf  u$ and $\delta \mathbf  u$ are the displacement trial and test functions, respectively. 
The traction, $ \mathbf  f _N$, is applied to the boundary $\Gamma_N$. The Cauchy stress tensor, $ \bm  \sigma$, is defined as $ \bm  \sigma = \bm C_{eff} \bm{\varepsilon}$ 
based on the assumption of
linear elasticity where  $\bm{\varepsilon}$ denotes the strain tensor and $\bm C_{eff}$ is the homogenized 
elasticity tensor obtained from Mori-Tanaka homogenization discussed in Section \ref{materialmodel}. Using Nitsche's unsymmetric formulation, 
the Dirichlet boundary condition is weakly imposed along the boundary $\Gamma_D$ as follows:

\begin{equation}
\begin{aligned}
\mathcal R^D = & -\int_{\Gamma^{D}} \delta \mathbf u \cdot \sigma ( \mathbf u ) \mathbf n _{\Gamma} d \Gamma \\
& +\int_{\Gamma^{D}} \mathbf \sigma (\delta  \mathbf u )  \mathbf n _{\Gamma} \cdot\left(  \mathbf u -  \mathbf u _D\right) d \Gamma \\
& + \gamma_{D} \int_{\Gamma^{D}} \delta \mathbf u \cdot\left( \mathbf u - \mathbf  u _D\right) d \Gamma ,
\end{aligned}
\end{equation}

\noindent where $\mathbf u _D$ is the prescribed displacement and $ \mathbf n _{\Gamma}$ vector 
denotes the normal pointing outward from the boundary $\Gamma_{D}$. 
The penalty factor depends on mesh size and is obtained as follows: 

\begin{equation}
\gamma_{D}=c_{D} \frac{E_{eff}}{h}, 
\end{equation}

\noindent where $c_{D}$ is a parameter that controls the accuracy of the Dirichlet enforcement, and $E_{eff}$ 
stands for effective Young's modulus of the FRC material, computed by homogenization.

The term \( \mathcal{R}^G \) denotes the facet-oriented ghost stabilization contribution.
 When the level set intersects background elements, it may create small subdomains, diminishing the support of the basis function. 
 This can lead to ill-conditioning of the linear system and potentially degrade the accuracy of the state variables and 
 their gradient approximations. 
 To mitigate this issue, this work adopts the approach of \cite{2020NoelEtAl} for stabilizing XIGA and introduces the following ghost stabilization term:

\begin{multline}
\mathcal{R}^G = \\
\sum_{i=1}^{N_{F}} 
\sum_{j \in J_{F,i}} 
\Bigg[ \sum_{k=1}^{p} \int_{F}  
\gamma_{G}^{\mathbf{a}}\,
\Big\llbracket \partial^{k}_{n}\, \delta \mathbf{u} \Big\rrbracket 
\cdot
\Big\llbracket \partial^{k}_{n}\, \mathbf{u} \Big\rrbracket\, 
d\Gamma\ \Bigg],
\label{eq:GhostStabilization}
\end{multline}

\noindent where \( J_{F,i} \) denotes the set of facets subject to the ghost penalty, \( N_{F} \) is the total number of such facets,
 and \( p \) represents the order of approximation. \( \gamma_{G}^{\mathbf{a}} \) is the ghost penalty parameter,
  $\partial^k_n$ is the $k^{th}$ order derivative in the normal direction of the facet, and $\Big\llbracket \ \Big\rrbracket$
 is the jump operator measuring the difference across the ghost facet.  This term penalizes discontinuities in the derivatives of the state variables across the entire element facet.

 In the optimization process, 
 the LSF may evolve to create isolated, topologically disconnected subregions with rigid body modes. 
 Unlike the approach in \cite{wei2010study} where the void phase is modeled as a soft material, we adopt the approach in \cite{geiss2018topology} and introduce a weak 
 elastic bedding for free-floating
 regions. Identification of these subregions is achieved by solving an auxiliary thermal-convection problem, where the temperature is prescribed at mechanical boundaries. 
 The convection effect ensures that all isolated regions assume ambient temperature which is used to activate the elastic bedding.
 The formulation for the elastic bedding residual is as follows:

\begin{equation}
\mathcal{R}^{S} = \int_{\Omega^m} \gamma_s \kappa_s (\delta \mathbf{u}) \mathbf{u} \, d\Omega,
\label{eq:SelectiveSpring}
\end{equation}

\noindent where \( \kappa_s = \frac{E_{eff}}{h^2} \), with \( E_{eff} \) representing effective Young's modulus and \( h \) the mesh size. 
The coefficient \( \gamma_s \) activates the elastic bedding and is defined as a smooth transition function of the auxiliary temperature, 
in the interval of $[0,T_{pre}]$ where $T_{pre}$ denotes the prescribed temperature.

\subsection{Mori-Tanaka homogenization}\label{materialmodel}

\begin{figure}[b]\centering
	\includegraphics[scale=1]{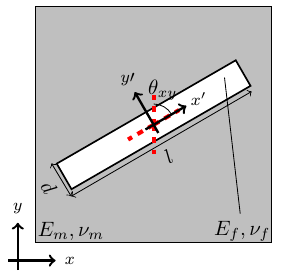}
	\caption{Matrix fiber layout for MT.}
	\label{fig_mt}
\end{figure}

To predict the physical response of the FRC, the Mori-Tanaka (MT) homogenization scheme \citep{mori1973average} is employed. 
MT is a mean-field homogenization technique that assumes uniform inhomogeneity in the matrix and builds on Eshelby's elasticity method. Figure \ref{fig_mt} 
demonstrates the idealized Representative Volume Element (RVE) for the MT homogenization scheme. It considers the properties of the individual phases, 
as well as the geometry and arrangement of the phases within the composite material. These properties include Young's 
modulus of the matrix and fiber, $E_m, E_f$, Poisson ratio of the fiber and matrix, $\nu_m, \nu_f$, volume fraction, $V_f$, and aspect ratio, $AR  = \frac{d}{l}$. 
The rotation angles, $\theta_{xy}$ and $\theta_z$, are used to construct the rotation tensor which transforms the stiffness tensor from the local coordinate system, aligned with the fiber orientation, to the global coordinate system.
The effective constitutive tensor ($\bm{C}_{eff}$) is expressed as a function of the fiber and
 matrix's material properties as well as the fiber orientation as:


\begin{align}
	\bm{C}_{eff} &= \nonumber \mathbf{Q}^T(\theta_{xy},\theta_z)\\
	&  \cdot \mathbf{C}_{MT}(E_m, E_f,\nu_m,\nu_f,V_f,AR) \cdot \mathbf{Q}(\theta_{xy},\theta_{z}),
	\label{eqFiberBspline}
\end{align}

\noindent where the rotation tensor is denoted by $\mathbf{Q}$.

%% file: Opt.tex
\section{Optimization framework}\label{optimzation}

The optimization problems addressed in this study are described by the following general formulation:

\begin{equation}
\begin{aligned}
\min_{\mathbf{s}} & \quad z(\mathbf{s}, \mathbf{u}(\mathbf{s})) \\
\text{s.t.} & \quad g_j(\mathbf{s}, \mathbf{u}(\mathbf{s})) \leq 0, \quad j=1,\ldots,N_g, \\
& \quad \underline{s} \leq s_i \leq \overline{s}, \quad i=1,\ldots,N_s,
\end{aligned}
\label{eq:optimization_problem_generic}
\end{equation}

\noindent where \( \mathbf{s} \) denotes the vector of optimization variables of dimension \( N_s \), 
constrained between lower bounds \( \underline{\mathbf{s}} \) and upper bounds \( \overline{\mathbf{s}} \). 
The state variable vector \( \mathbf{u}(\mathbf{s}) \) denotes the structural displacements. The objective function is represented by \( z(\mathbf{s}, \mathbf{u}(\mathbf{s})) \), 
with \( g_j(\mathbf{s}, \mathbf{u}(\mathbf{s})) \) as the constraint functions. The focus of this paper is on the minimization of compliance, 
augmented by the regularization terms. The objective function is defined as follows:

\begin{equation}
\begin{aligned}
z(\mathbf{s}, \mathbf{u}(\mathbf{s})) &= 
w_f \frac{\mathcal{F}(\mathbf{s}, \mathbf{u}(\mathbf{s}))}{\mathcal{F}^0(\mathbf{s}, \mathbf{u}(\mathbf{s}))} 
+ w_p \frac{\mathcal{P}_{p}(\mathbf{s})}{\mathcal{P}^0_{p}(\mathbf{s})}
+ w_g \frac{\mathcal{P}_{g}(\mathbf{s})}{\mathcal{P}^0_{g}(\mathbf{s})} \\
&+ w_{\text{par}} \frac{\mathcal{P}_{\text{par}}(\mathbf{s})}{\mathcal{P}^0_{\text{par}}(\mathbf{s})}
+ w_{\text{Lcur}} \frac{\mathcal{P}_{\text{Lcur}}(\mathbf{s})}{\mathcal{P}^0_{\text{Lcur}}(\mathbf{s})} \\
&+ w_{\text{Gcur}} \frac{\mathcal{P}_{\text{Gcur}}(\mathbf{s})}{\mathcal{P}^0_{\text{Gcur}}(\mathbf{s})},
\end{aligned}
\label{eq:objective_function}
\end{equation}

\noindent where \( \mathcal{F}(\mathbf{s}, \mathbf{u}(\mathbf{s})) \) denotes the strain energy. 
The optimization problem is regularized by perimeter and level set gradient penalties, \( \mathcal{P}_{p} \) and \( \mathcal{P}_{g} \), 
with weights \( w_p \) and \( w_g \), respectively. The alignment and curvature of the fiber orientation fields are controlled by \( \mathcal{P}_{\text{par}} \), 
\( \mathcal{P}_{\text{Lcur}} \) and \( \mathcal{P}_{\text{Gcur}} \), representing the parallel misalignment penalty, local curvature penalty and
global curvature penalty, with \( w_{\text{par}} \), \( w_{\text{Lcur}} \) and \( w_{\text{Gcur}} \) being the associated weights.
All terms in the objective function are normalized by reference values denoted by the superscript \( 0 \). Subsequent subsections will discuss these penalties in greater detail. 

\subsection{Regularization of the Level Set}\label{RGLS}

\cite{haber1996new} introduced a constraint on the perimeter of the solid-void interface to discourage 
irregular geometrical features. 
Following this approach, a perimeter penalty term is incorporated into the objective function as follows:

\begin{equation}
\mathcal{P}_{p} = \frac{1}{\mathcal{P}_0} \int_{\Gamma_I}\, d\Gamma,
\label{eq:perimeter_penalty}
\end{equation}

\noindent where \( \mathcal{P}_0 \) denotes the initial perimeter value. 

In addition, a regularization term, \( \mathcal{P}_{g} \), is employed to stabilize the gradient of LSF in proximity to the interface.
 This term promotes uniform spatial gradients of the LSF in the vicinity of zero iso-contour. As the LSF evolves throughout the optimization process, convergence to upper and lower 
 bounds is encouraged away from the 
 interface. This approach prevents the LSF from becoming excessively flat or steep which may
 cause oscillations in the optimization process. 
 The regularization approach outlined by \cite{geiss2019regularization} introduces a smoothed target field with uniform gradients along the interface. Deviations from this target field in terms of function values and gradients incur a penalty in the objective function:

 \begin{equation}
\mathcal{P}_{r} = \frac{\int_{\Omega} w_{\phi} \left( \phi-\tilde{\phi} \right)^2 d\Omega}{\int_{\Omega} \phi_{Bnd}^2\, d\Omega}
+ \frac{\int_{\Omega} w_{\nabla \phi} \left| \nabla \phi - \nabla\tilde{\phi} \right|^2 d\Omega}{\int_{\Omega} d\Omega}.
\label{eq:regularization_penalty}
\end{equation}

Here \( \tilde{\phi} \) denotes the corresponding smoothed target field. The penalty weights \( w_{\phi} \) and \( w_{\nabla \phi} \) are responsible for penalizing deviations in the level set values and its gradient. These weights are adjusted based on the proximity to the interface.
The level set lower bound is denoted by \( \phi_{Bnd} \), and the target field \( \tilde{\phi} \) is constructed as a truncated signed distance field through a sigmoidal function:

\begin{equation}
\tilde{\phi} = \left( \frac{2}{1+\exp({-2\, \phi_{D}} / { \phi_{Bnd}})} -1 \right)\,  \phi_{Bnd},
\label{eq:target_field}
\end{equation}

\noindent with \( \phi_{D} \) as the signed distance field, computed by the "heat method" of \cite{crane2017heat}.

\subsection{Parallel misalignment penalty}\label{par_pen}
This section introduces a key contribution of the paper: 
the misalignment penalty term, which facilitates the generation of locally parallel fiber paths. 
Such parallelism minimizes the post-processing required to transform the fiber orientation field into continuous fiber paths, 
promoting a more seamless transition.

Geometrically, parallel curves can be defined in several different ways. Two curves are parallel if the distance between them is constant at 
every point along their length. 
Alternatively, two curves within the same plane can be considered parallel if they have identical 
tangent vectors at corresponding points along their respective lengths.
This study takes advantage of the second definition, and curves are defined as parallel 
when the tangent does not change in the normal direction of the curve. 

\begin{figure}[!h]\centering
	\includegraphics[scale=1]{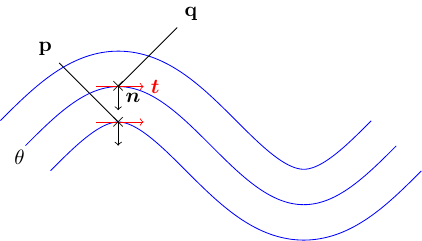}
	\caption{Parallel fiber path.}
	\label{fig_streamline}
\end{figure}

Consider Fig.~\ref{fig_streamline} which presents three geometrically parallel curves. 
Traversing from point $q$ in the normal direction, $\bm{n}$, leads to point $p$. The tangent vectors in both points are identical.
Consequently, the tangent vector $\bm{t}$ remains invariant in the direction of the normal vector $\bm{n}$. 
In mathematical terms, the directional derivative of $\bm{t}$ with respect to $\bm{n}$ equates to zero, i.e.:

\begin{equation}
 \nabla \bm{t} \cdot \bm{n} = \bm{0}.
 \label{eqParallel}
\end{equation}

\begin{figure}[b!]\centering
	\includegraphics[scale=2.0]{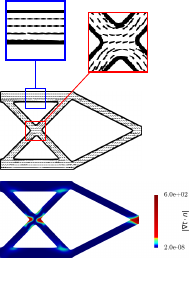}
	\caption{Misalignment penalty visualization for a specific fiber orientation field.}
	\label{fig_parallel}
\end{figure}

Figure~\ref{fig_parallel} illustrates the fiber orientation field's layout at the top, with each line depicting a fiber tangent as defined by Equation (\ref{eqTangentVector}).
The corresponding \( \nabla \bm{t} \cdot \bm{n} \) values for the fiber orientation field are shown in the contour plot at the bottom. These values are zero in regions with parallel fiber paths and non-zero where fibers intersect.

For a 2D problem, based on the parameterization of the fiber orientation in Equation (\ref{eqTangentVector}), 
there is a direct relationship between the tangent vector and the normal vector as well as the angle value. As the tangent and normal are constructed from $\theta_{xy}$, 
for parallel curves, one can infer that the change in $\theta_{xy}$ in the direction of $\mathbf n$ is zero, i.e.:

\begin{equation}
\nabla \bm{t} \cdot \bm{n}  \equiv  \nabla \bm{n} \cdot \bm{n}  \equiv  \nabla \theta_{xy} \cdot \bm{n} =  \bm{0}.
\label{eqEquivalent}
\end{equation}

In 3D, with two angles defining the tangent, and the normal vector being non-unique, 
the condition \( \nabla \bm{t} \cdot \bm{n} = 0 \) still applies, but \( \bm{n} \) could be any vector with its normal vector being \( \bm{t} \). 
A vector in this plane can be decomposed into two linearly independent vectors \( \bm{n}_1 \) and \( \bm{n}_2 \) with arbitrary coefficients \( a \) and \( b \) as:
$ \bm{n} = a \bm{n}_1 + b \bm{n}_2.$
Subsequently, Equation (\ref{eqParallel}) can be reformulated as: 

\begin{equation}
\nabla \bm{t} \cdot \bm{n} = a \left(\nabla \bm{t} \cdot \bm{n_1} \right) + b \left( \nabla \bm{t} \cdot \bm{n_2} \right) = \bm{0},
\end{equation}

By following the right-hand rule for \( \bm{n}_1 \), \( \bm{n}_2 \), and \( \bm{t} \), with \( \bm{n}_1 \) arbitrarily chosen to be orthogonal to 
\( \bm{t} \), and \( \bm{n}_2 \) obtained from their cross product, the parallel condition can be satisfied.

The parallel misalignment penalties for 2D and 3D structures are defined as:
\begin{equation}\label{eq:parallel_penalty2d3d}
	\begin{aligned}
		\mathcal{P}_{\text{par}}^{\text{2D}} & = \int_\Omega \| \nabla \bm{t} \cdot \bm{n} \|_2^2  dV, \\
		\mathcal{P}_{\text{par}}^{\text{3D}} & = \int_\Omega \| \nabla \bm{t} \cdot \bm{n_1} \|_2^2  dV + \int_\Omega \| \nabla \bm{t} \cdot \bm{n_2} \|_2^2  dV.
	\end{aligned}
\end{equation}

A simple optimization problem is presented to illustrate the impact of the parallel misalignment penalty. 
Here, the fiber orientation is the sole design parameter field, and the sole objective is to minimize parallel misalignment 
penalty \( \mathcal{P}_{\text{par}} \). The algorithmic setting uses the parameters outlined in Section \ref{sectionExamples}.

\begin{figure}[t]
	\centering
	\includegraphics[scale=1]{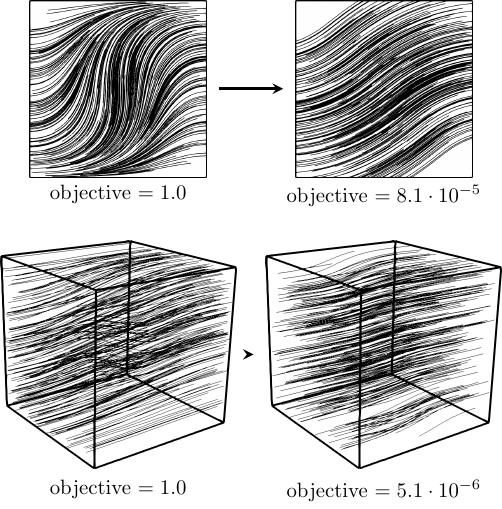}
	\caption{Minimization of the parallel misalignment penalty for 2D and 3D: initial layout (left) and optimized layout (right).}
	\label{fig_par_pen}
\end{figure}

Figure \ref{fig_par_pen} shows the initial and optimized fiber layouts for both 2D and 3D cases. For visualization purposes, the fiber paths are generated by equally spaced streamlines.
 The initial layout is generated by the function 
\( \theta_{xy}(x,y,z) = \theta_z(x,y,z) = \sin(\pi x)\sin(\pi y)\sin(\pi z) \) for 3D, and 
\( \theta_{xy}(x,y) = \sin(\pi x)\sin(\pi y) \) for 2D. Visual inspection suggests that the parallel penalty effectively aligns the fiber paths.
The normalized objective value (parallel misalignment penalty), relative to the initial design, 
reduces from $1.0$ to $8.1\cdot 10^{-5}$ for the 2D case, and from $1.0$ to $5.1\cdot 10^{-6}$ in the 3D case.

\subsection{Curvature penalty}

The misalignment penalty introduced in Equation (\ref{eq:parallel_penalty2d3d}) 
promotes the parallel fiber paths but does not control the curvature of the fibers.
 To avoid excessive bending of the fibers, 
a penalty term is introduced that limits the curvature of the fiber path.
The curvature 
\(\kappa(l)\) for a curve parameterized by arc length \(l\) is:

\begin{equation}
\kappa(l) = \frac{\lVert \bm{r}'(l) \times \bm{r}''(l) \rVert}{\lVert \bm{r}'(l) \rVert ^3 },
\end{equation}

\noindent where the position vector is denoted by \(\bm{r}(l)\), with \(\bm{r}'(l)\) and \(\bm{r}''(l)\) 
representing its first and second-order derivatives with respect to the arc length, respectively. 
The second-order derivative of the position vector, derived using the chain rule, is 
\(\bm{r}''(l) = \frac{D \bm{r}'(l)}{D \bm{x}} \cdot \bm{r}\). For fiber paths generated from fiber orientation field(s),
 the first order derivative of the position vector \(\bm{r}'\) is equal the tangent vector \(\bm{t}\), 
 given in Equation (\ref{eqTangentVector}), i.e., \(\bm{r}' =\bm{t}\). Using the chain rule, the second 
 derivative of the position vector is computed as: \(\bm{r}''  = \bm{t}' =  \frac{D \bm{t}}{D \bm{x}} \cdot \bm{t}\).

Given the unit magnitude of the tangent vector, 
$ \lVert \mathbf{t} \rVert = 1,$
the curvature of the fiber path is subsequently calculated as:

\begin{equation}
\kappa^{2D,3D} = \lVert \bm{M} \bm{t} \times \bm{t} \rVert,
\end{equation}

\noindent where the matrix \(\bm M\) represents
the spatial derivative of the tangent vector field and is defined as \( \bm M = \frac{D \bm{t}}{D \bm{x}}\). 
For 2D problems,
\( \bm M\) is a function of \(\theta_{xy}(x,y)\) while for 3D problems it depends on
both \(\theta_{xy}(x,y,z)\) and \(\theta_{z}(x,y,z)\).

\begin{figure}[b]
	\centering
	\includegraphics[scale=2.0]{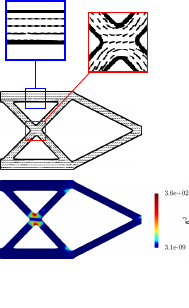}
	\caption{Curvature penalty visualization for a specific fiber orientation field.}
	\label{fig_par_cur}
\end{figure}

In Fig.~\ref{fig_par_cur}, the upper section displays the configuration of the fiber orientation field, 
where individual lines denote fiber tangents as defined in Equation (\ref{eqTangentVector}). 
The corresponding \( \kappa^2 \) 
values for the fiber orientation field are shown at the bottom.
These values are zero in regions with straight fibers and non-zero where fibers exhibit bending, 
as illustrated by the blue and red magnified insets in the figure, respectively.

To ensure local curvature does not exceed a maximum manufacturing limit defined by \(\kappa_{\text{max}}\), 
the following penalty term is incorporated into the objective function:

\begin{equation}\label{eq:curvature_penalty}
	\mathcal{P}_{\text{Lcur}}^{\text{2D,3D}} = \int_\Omega \left( \left( \kappa^2 - \kappa_{\text{max}}^2 \right)^{+} \right)^2 dV, 
\end{equation}

\noindent where $( \ . \ )^+  = max(0, \ . \ )$. The penalty term agglomerates the local point-wise curvature constraint into a global 
constraint which is applied as a penalty term in the objective function.

While the local curvature penalty aims to lessen sharp turns in fiber paths, 
a global curvature penalty can be formulated to enhance overall smoothness, 
discouraging wavy fibers. The global curvature penalty is defined as follows:

\begin{equation}\label{eq:gcurvature_penalty}
	\mathcal{P}_{\text{Gcur}}^{\text{2D,3D}} = \int_\Omega \kappa^2  dV. 
\end{equation}

In 2D problems, the curvature is simplified by substituting the tangent vector from Equation 
(\ref{eqTangentVector}) into Equation (\ref{eq:curvature_penalty}), yielding \(\nabla \theta_{xy} \cdot \bm{t} \). 
This represents the rate of change in the fiber angle along the tangent, capturing the fiber path's curvature. 
With this simplification, the following relation holds:

\begin{equation}
\| \nabla \theta_{xy} \|^2 = \| \nabla \theta_{xy} \cdot \bm{t} \|^2 + \| \nabla \theta_{xy} \cdot \bm{n} \| ^2.
\end{equation}

This demonstrates that combining misalignment and
curvature penalties with equal weights effectively penalizes the first-order derivative of the fiber orientation field.
Adjusting the weights or using anisotropic penalties can
produce outcomes not achievable through penalization
of the first derivative alone. This effect is illustrated in
the previous section, where the parallel misalignment penalty was minimized, and
further explored in the rest of the section.

To demonstrate the curvature penalty's effect on fiber layout smoothness, we consider an optimization problem where the algorithmic setting uses the parameters outlined in Section \ref{sectionExamples}. 
Here, the objective is to minimize the global curvature penalty, with the fiber orientation being the only design variable. 

\begin{figure}[!ht]
	\centering
	\includegraphics[scale=1]{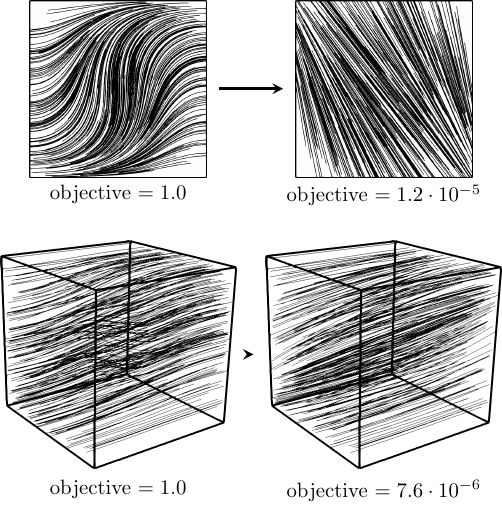}
	\caption{Minimization of the curvature penalty for 2D and 3D: initial layout (left) and optimized layout (right).}
	\label{fig_laplace_pen}
\end{figure}

The initial and optimized fiber layouts are illustrated in Fig.~\ref{fig_laplace_pen}. 
The optimization results in straight fiber paths with zero curvature,
as evident by the figure, 
with the normalized objective value with respect to the initial design reducing from \( 1.0 \) to \( 1.2\cdot 10^{-5} \) in 2D
case and from  \( 1.0 \) to \( 7.6\cdot 10^{-6} \) in 3D case.
It is important to note that while the optimized fiber path is straight, it is not parallel due to the exclusion of the parallel misalignment penalty in the optimization problem.

%% file: Examples.tex
\section{Numerical optimization examples}\label{sectionExamples}

The subsequent sections study the interplay between geometry and fiber orientation field and the efficacy 
of the proposed methodology for achieving parallel 
and smooth fiber orientation configurations. 
Initially, a 2D optimization problem is presented with fiber orientation as the sole design parameter to 
isolate the impact of misalignment and curvature penalties. This example further explores the effect of 
the polynomial order of the fiber orientation parameterization. 
Subsequently, an optimization problem considering both structural shape and fiber orientation in 2D is considered, 
investigating the influence of fiber orientation parameterization in terms of the refinement level. 
Lastly, two three-dimensional cases are examined to optimize geometry and fiber orientation simultaneously, 
analyzing the effect of misalignment and curvature penalties in 3D settings.

\begin{table}[h!]
	\centering
\begin{tabular}{|c | c|} 
	\hline
	Parameter & Value  \\ 
	\hline
	$E_m$  & 1.03   \\  \hline
	$E_f$  & 1.02 $\cdot$ $10^{3}$  \\ \hline
	$\nu_m$  & 0.4 \\ \hline
	$\nu_f$  & 0.4 \\ \hline
	$AR$  & 10.0  \\ \hline
	$V_f$  & 0.1  \\
	\hline
\end{tabular}
\caption{Matrix and fiber material properties.}
\label{tablematprop}
\end{table}

Table~\ref{tablematprop} provides the material properties used in the numerical examples. This includes the matrix and fiber Young's moduli, 
their respective Poisson's ratios, fiber aspect ratio and volume fraction. To amplify the anisotropic behavior of the fibers and their influence on the structural response, a combination of high stiffness ratio and low volume fraction is selected. Units are omitted as they are consistent throughout the study.

In this study, our computational domains are defined on rectangular grids, 
utilizing uniform B-spline meshes for discretization. Mesh refinement is conducted by 
recursively subdividing each element into four elements in 2D and into eight elements in 3D until the specified level of refinement is attained.

The state variable field, i.e. the displacement field is discretized using bi-linear B-splines with quadrilateral elements for 2D, 
and tri-linear B-splines with hexahedral elements in 3D. The level set field is discretized with bi-quadratic (2D) and tri-quadric (3D) B-splines on a mesh that is twice as coarse as 
the state variable mesh. 
The fiber orientation field is discretized on a mesh that is either two or four times coarser than the state variable mesh, 
utilizing linear, quadratic and cubic B-splines. 
As discussed in \cite{2020NoelEtAl}, defining the design variables fields on coarser but higher-order B-spline meshes has a smoothing effect, 
suppressing spurious geometric features and oscillating fiber orientation fields.

For the evaluation of the governing equations, the weak Dirichlet boundary condition penalty, $\gamma_D$, is assigned a value of $10.0$, while the ghost penalty, 
$\gamma_G $, is chosen to be $0.01$. The discretized governing equations and adjoint sensitivity equations for 2D problems are solved using the direct solver 
PARDISO \citep{schenk2001pardiso}. 
In 3D problems, the linear systems are solved using a Generalized Minimal RESidual (GMRES) method, with an ILUT (dual threshold incomplete LU factorization) preconditioner. 
The convergence tolerance for the GMRES algorithm, i.e. the required drop in relative preconditioned residual, is set to $1\cdot 10^{-9}$.

The optimization problem is solved using a gradient-based algorithm, 
the Globally Convergent Method of Moving Asymptotes (GCMMA) referenced in \cite{Svanberg2002}, 
utilizing two inner iterations. In cases where the local curvature penalty is applied and
 both geometry and fiber orientation are treated as design variables, it was observed that the use of inner iterations
 resulted in designs that stagnate in a local minimum. Consequently, for such cases, 
 the optimization initially progresses without inner iterations for a specified number of steps, after which inner iterations are activated.
 The optimization is considered to have converged when the relative change in objective function values between two successive iterations falls below \(1 \cdot 10^{-5}\).
In the GCMMA algorithm, the parameters for adapting the initial, shrinking, and expanding asymptotes are set to $0.5$, $0.7$, and $1.2$, respectively.
The design sensitivities are computed by the adjoint method, see \cite{sharma2017shape} for more details on the adjoint method for XFEM problems.

Fiber orientation design variable fields, $\theta_{xy}$ and $\theta_z$, are bounded by the box-constraints $ \left[ -3\pi, 3\pi \right] $ for all the examples. The value for maximum feasible curvature, \(\kappa_{\text{max}}\) is selected to demonstrate the effectiveness of the local curvature penalty in reducing fiber path curvature. 
This parameter is set to a value much lower than the curvature values observed when only the parallel misalignment penalty is applied. 

We employ continuation strategies to regulate the influence of the misalignment and curvature penalties in the objective function.
By renormalizing the terms in the objective function, the contribution of the penalty terms is kept minimal, ensuring
that the strain energy is the predominant term in the objective function. 
The strain energy is initially normalized against its initial design value, 
whereas the penalty terms corresponding to the fiber orientation field are normalized based on their values at the 
tenth optimization iteration. 
The fiber paths are initially assumed to be horizontal or vertical across the domain, unless stated otherwise, leading to initial penalty values being zero. 
Consequently, the penalty terms are activated after ten optimization iterations, 
utilizing their corresponding values at this stage as normalization factors.
Renormalization occurs at every 10th outer iteration step. 
The local curvature penalty weight is chosen to be $0.1$ and every 20 optimization steps this weight is multiplied by 4.0. 
The weights assigned to the misalignment penalty and the global curvature penalty are \(0.05\) and \(0.01\), respectively.   

\subsection{Post-processing for fiber paths}
In our methodology, the fiber orientation field is directly used as a design variable in the optimization process.
This parametrization does not immediately
provide a geometrical description of the continuous fiber paths. Therefore,
a post-processing step is needed to convert the fiber orientation
field into a geometric description of the fiber path while preserving the volume fraction. We employ the methodology in
 \cite{boddeti2020optimal} to determine continuous fiber paths. This process uses the
stripe patterns algorithm by \cite{knoppel2015stripe}, which
positions evenly spaced and parallel stripes along a specified
vector field on a manifold, with singularities introduced
as needed to maintain parallelism and even
spacing. Alternatively to visualizing the fiber orientation fields via stripe patterns, 
we visualize the raw optimization results via the streamline functionality in Paraview \citep{Moreland}, 
where tangents to the fiber fields are illustrated.
Thus, the length of the tangent vector in the streamline visualization does not have any physical meaning.

\subsection{Short plate under shear}\label{example1}

\begin{figure}[t]\centering
	\typeout{>>>>>>\the\linewidth}
	\includegraphics[width=\linewidth]{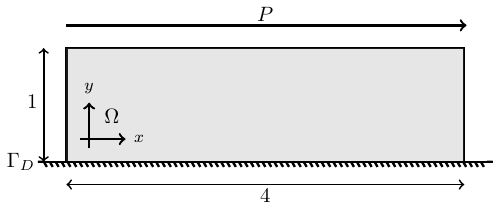}
	\caption{Short plate configuration.}
	\label{fig_short_cantilver_beam_44}
\end{figure}

This example investigates how various combinations of misalignment and curvature penalties affect the optimized 
fiber layout. 
We explore different scenarios employing linear, quadratic, and cubic polynomial orders for the
fiber orientation B-spline discretization. 
The sole design parameter in this problem is the fiber orientation, while the structural shape remains unchanged. 
This example minimizes the compliance of the short plate shown in Fig.~\ref{fig_short_cantilver_beam_44} 
under plane stress conditions.
The plate is fixed at the bottom and a load ($P = 0.01$) is applied at the top. 
The short geometry of the plate induces shear-dominated behavior in its structural response.
 
With the structural shape fixed, the formulation of the optimization problem in 
Equation (\ref{eq:objective_function}) 
is simplified to an unconstrained optimization problem as follows:

\begin{equation}
\begin{aligned}
\underset{\theta_{xy}}{min} \Pi  = & w_{f} \frac{\strainEnergy}{\strainEnergy^0}  +
 \frac{w_{\text{par}}}{\mathcal{P}_{\text{par}}^{0} } \mathcal{P}_{\text{par}}  + \\ 
&  \frac{w_{\text{Lcur}}}{\mathcal{P}_{\text{Lcur}}^{0} }  \mathcal{P}_{\text{Lcur}} + \frac{w_{\text{Gcur}}}{\mathcal{P}_{\text{Gcur}}^{0} }  \mathcal{P}_{\text{Gcur}}.
\end{aligned}
\label{Eq_example1}
\end{equation}

The first term in the Equation (\ref{Eq_example1}) captures the strain energy of the system 
and $\strainEnergy^{0}$ denotes the initial strain energy.
The rest of the terms correspond to the misalignment and curvature penalties 
defined in Equation (\ref{eq:parallel_penalty2d3d}), Equation (\ref{eq:curvature_penalty}) and Equation (\ref{eq:gcurvature_penalty}). 
All terms are normalized by reference values indicated by the superscript \(0\). The maximum allowable curvature is set to $\kappa_{\text{max}} = 1.0$.

To analyze the effects of the misalignment and curvature penalties, four different penalization cases are considered: 
\begin{itemize}
    \item Case NP: No penalties applied - $w_{\text{par}} = 0.0$, $w_{\text{Lcur}} = 0.0$, $w_{\text{Gcur}} = 0.0$.
    \item Case MP: Only the misalignment penalty is applied - $w_{\text{par}} = 0.05$, $w_{\text{Lcur}} = 0.0$, $w_{\text{Gcur}} = 0.0$.
    \item Case MLCP: Both misalignment and local curvature penalties are applied - $w_{\text{par}} = 0.05$, $w_{\text{Lcur}} \ne 0.0$, $w_{\text{Gcur}} = 0.0$.
    \item Case MCP: All penalties, including misalignment, local curvature, and global curvature, are applied - $w_{\text{par}} = 0.05$, $w_{\text{Lcur}} \ne 0.0$, $w_{\text{Gcur}} = 0.01$.
\end{itemize}

The fiber orientation field is discretized on linear, quadratic or 
cubic B-spline meshes. To obtain a filtering effect, 
these B-spline meshes are coarser than the linear B-spline meshes discretizing the state variable field. 
In this example, the fiber orientation mesh is selected to be two times coarser than the state variable mesh, 
which uses a resolution of \(64 \times 256\) elements for all testing configurations. 
Consequently, the fiber orientation mesh resolution is set to \(16 \times 64\). The initial fiber orientation is set
to be vertical across the domain.  

\begin{figure}[b]\centering
	\includegraphics[width=0.4\textwidth]{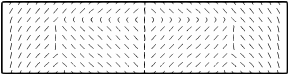}
	\caption{Streamline visualization of the optimized fiber path for quadratic B-splines discretization of the fiber orientation field for Case NP.}
	\label{fig_short_cantilver_beam_ref}
\end{figure}

Figure \ref{fig_short_cantilver_beam_ref} 
presents the optimized fiber path for a quadratic B-spline discretization 
under Case NP, where no penalty is applied, and Table \ref{EX1_ref_plate} displays the absolute values for the penalty terms, strain energy, and maximum curvature.
This configuration acts as the reference for subsequent analyses.

\begin{table}[ht]
    \centering
    \begin{tabular}{|c|c|}
    \hline
    Parameter & Value \\\hline
    $\strainEnergy$ & $3.925\cdot 10^{-4}$ \\ \hline
    max($\kappa$) & $37.09$ \\ \hline
    $\mathcal{P}_{par}$ & $110.513$ \\ \hline
    $\mathcal{P}_{Lcur}$ & $7.589\cdot 10^{4}$ \\ \hline
    $\mathcal{P}_{Gcur}$ & $1.199\cdot 10^{2}$ \\ \hline
    \end{tabular}
    \caption{Absolute values of objective components and maximum curvature for quadratic B-splines discretization of the fiber orientation field
	under Case NP.}
    \label{EX1_ref_plate}
\end{table}

\begin{figure*}[!h]\centering
	\includegraphics[scale=1.0]{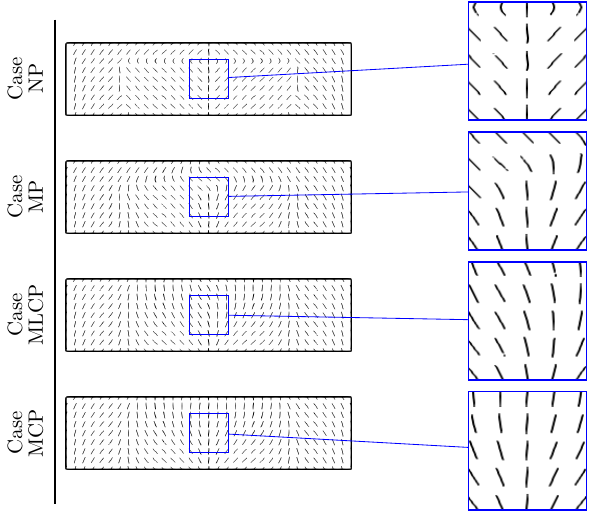}
	\caption{Streamline visualization of the optimized fiber layout of the short plate for different penalization cases with quadratic B-splines discretization of the fiber orientation field.}
	\label{fig_plate_with_hole_results}
\end{figure*}

\begin{figure*}[!h]\centering
	\includegraphics[scale=1.0]{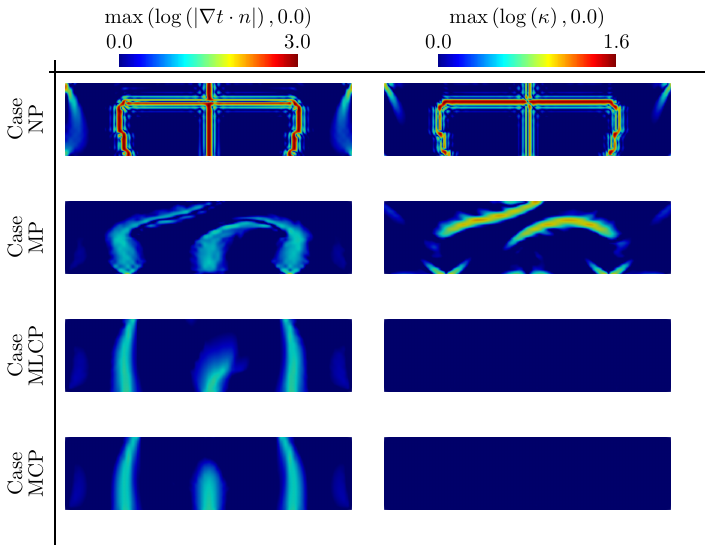}
	\caption{Misalignment penalty and curvature contour plots for different penalization cases for a quadratic B-spline fiber orientation field.}
	\label{fig_plate_with_hole_resultscontour}
\end{figure*}

Figure \ref{fig_plate_with_hole_results} shows the optimized fiber layouts obtained using 
quadratic B-splines for all cases, incorporating different penalty combinations. 
The corresponding strain energy value and penalty values, $\mathcal{P}_{\text{par}}$, $\mathcal{P}_{\text{Lcur}}$ and 
$\mathcal{P}_{\text{Gcur}}$, all normalized against the reference case and maximum curvature are displayed in Table \ref{EX1_Table_Order}.
Figure \ref{fig_plate_with_hole_resultscontour} presents the local value of the misalignment penalty and curvature values across 
different penalization cases, employing $max\left(\log_{10}(\cdot), 0.0\right)$ for better visualization. 
The $0.0$ threshold in the visualization function for curvature corresponds to $\log_{10}(\kappa_{max}) = \log_{10}(1.0) = 0.0$, 
whereas the threshold for the misalignment penalty is selected purely for visualization purposes.

\begin{table*}[!h]
	\centering
	\begin{tabular}{|c|c|c|c|c|c|}
	\hline
	Case & $\strainEnergy$ & max($\kappa$) & $\mathcal{P}_{par}$ & $\mathcal{P}_{Lcur}$ & $\mathcal{P}_{Gcur}$ \\
	\hline
	NP & 1.000 & 37.09 & 1.000 & $1.000$ & $1.000$ \\ \hline
	MP & 1.051 & 12.60 & 0.071 & $1.740\cdot 10^{-2}$ & $1.899\cdot 10^{-1}$ \\ \hline
	MLCP & 1.083 & 1.05 & 0.076 & $3.270\cdot 10^{-10}$ & $2.186\cdot 10^{-2}$ \\ \hline
	MCP & 1.092 & 1.03 & 0.068 & $3.307\cdot 10^{-12}$ & $1.651\cdot 10^{-2}$ \\ 
	\hline
	\end{tabular}
	\caption{Objective components, normalized against the reference case, and maximum curvature for different penalization cases with a quadratic B-spline fiber orientation field.}
	\label{EX1_Table_Order}
\end{table*}

It is observed that the misalignment penalty effectively aligns the fibers, 
reducing intersections of the fiber paths, 
as highlighted in the magnified insets in Fig.~\ref{fig_plate_with_hole_results}. 
The comparison of the parallel misalignment penalty in Table \ref{EX1_Table_Order}, $\mathcal{P}_\text{par}$, 
shows a substantial decrease in cases where the penalty is applied, 
indicating a higher degree of parallelism in fiber paths. This observation is further substantiated by the contour plots of 
the misalignment penalty in Fig.~\ref{fig_plate_with_hole_resultscontour},
which not only show decreased penalty values but also a reduction in regions of fiber intersection.

Incorporating local curvature penalty reduces the bending in the fiber paths, as illustrated in the magnified insets 
in Fig.~\ref{fig_plate_with_hole_results}
 , and decreases 
the curvature values close to the maximum allowable curvature value of 1.0 as demonstrated in Table \ref{EX1_Table_Order}.
Contour plots of the curvature values for Case MLCP and MCP in Fig.~\ref{fig_plate_with_hole_resultscontour}
demonstrate that curvature values, processed with the function $max\left(\log_{10}(\cdot), 0.0\right)$, are close to zero.
This implies that curvature values across the computational domain fall below or are very close to 1.0, 
indicating that the local curvature constraint is tightly satisfied.

The global curvature penalty further makes the fiber paths more uniform and straighter, as evident from the 
magnified insets of the fiber path in Fig.~\ref{fig_plate_with_hole_results} and the contour plots 
in Fig.~\ref{fig_plate_with_hole_resultscontour}. 
When comparing the curvature penalty values, $\mathcal{P}_{Lcur}$ and $\mathcal{P}_{Gcur}$, presented in Table 
\ref{EX1_Table_Order} 
it is evident that the two curvature penalties
are interconnected; applying one results in a decrease
in the other.

The addition of either misalignment or curvature penalties results in 
some non-symmetry in the optimized fiber layout, as can be seen in the magnified insets in 
Fig.~\ref{fig_plate_with_hole_results}. This non-symmetry arises because the penalty terms are not symmetric, 
resulting in non-symmetric gradients. 
For instance, for a symmetric function \(\theta(x,y) = x^2 + y^2\) in the computational domain \(\Omega = [-1,1] \times [-1,1]\),
the misalignment penalty \linebreak[4] \(\nabla \theta \cdot \bm{n} = 2x \cdot \cos\left(x^{2}+y^{2}\right) + 2y \cdot 
\sin\left(x^{2}+y^{2}\right)\) 
creates a non-symmetric objective function, which leads to non-systematic gradients.

Comparing strain energy values across different cases, presented in Table \ref{EX1_Table_Order} reveals that the addition 
of the misalignment and curvature penalties to the objective leads to increased strain energy. In a pure shear problem, 
optimal compliance corresponds to fiber orientations at \(\frac{\pi}{4}\) and \(\frac{3\pi}{4}\), 
aligning with the principal directions of stress. The constitutive tensor's shear components are identical at these angles. 
Given the shear-dominated nature of this problem, 
angles close to these values are observed in the optimized fiber layout for all cases with a 
continuous transition between these angles. This can be seen in the magnified insets in Fig.~\ref{fig_plate_with_hole_results}.

\begin{figure*}[!h]\centering
	\includegraphics[scale=1.0]{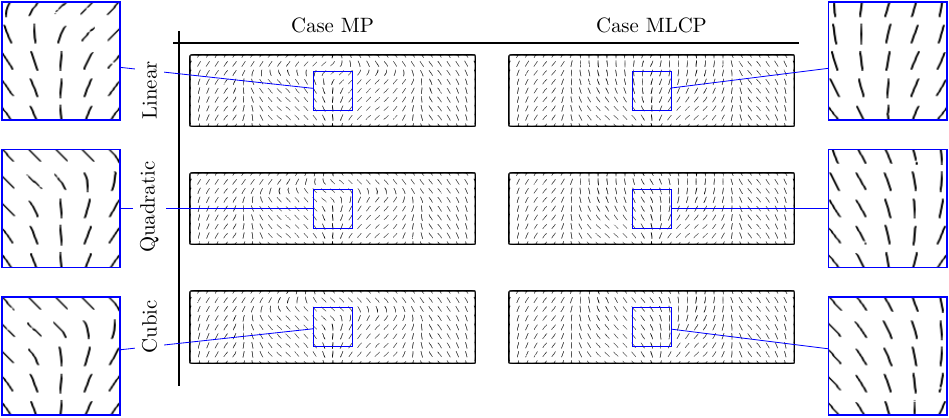}
	\caption{Streamline visualization of the optimized fiber layout of the short plate for different polynomial order of the fiber orientation field.}
	\label{fig_short_cantilver_beam_cubic}
\end{figure*}

\begin{table*}[ht]
    \centering
    \begin{tabular}{|c||c|c|c||c|c|c||}
    \hline
     & \multicolumn{3}{c||}{Case MP} & \multicolumn{3}{c||}{Case MLCP} \\
    \hline
     & Linear & Quadratic & Cubic & Linear & Quadratic & Cubic \\
    \hline
    $\strainEnergy$ & 1.053 & 1.051 & 1.051 & 1.092 & 1.083 & 1.079 \\
    \hline
    max($\kappa$) & 12.89 & 12.60 & 12.60 & 1.16 & 1.05 & 1.03 \\
    \hline
    $\mathcal{P}_{\text{par}}$ & 0.083 & 0.071 & 0.068 & 0.084 & 0.076 & 0.058 \\
    \hline
    $\mathcal{P}_{\text{Lcur}}$ & $2.191\cdot 10^{-2}$ & $1.740\cdot 10^{-2}$ & $1.174\cdot 10^{-2}$ & $1.077 \cdot 10^{-9}$ & $3.270 \cdot 10^{-10}$ & $1.545 \cdot 10^{-10}$ \\
    \hline
    $\mathcal{P}_{\text{Gcur}}$ & $1.953\cdot 10^{-1}$ & $1.899\cdot 10^{-1}$ & $1.782\cdot 10^{-1}$ & $2.074\cdot 10^{-2}$ & $2.186\cdot 10^{-2}$ & $1.682\cdot 10^{-2}$ \\
    \hline
    \end{tabular}
	\caption{Objective components, normalized against the reference case, and maximum curvature for different penalization cases and discretization orders.}
    \label{cases_by_refinement}
\end{table*}

To examine the impact of the polynomial order, we consider Cases MP and MLCP and conduct simulations with varying discretization orders for the fiber orientation field. 
Fig.~\ref{fig_short_cantilver_beam_cubic} shows the optimized fiber layouts using linear, quadratic, and cubic B-splines 
with a mesh size of $16 \times 64$. The corresponding objective components, normalized against the reference case and maximum curvature, 
are presented in Table \ref{cases_by_refinement}.

No visual distinctions are noted between the fiber paths in the quadratic and cubic cases, 
whereas the linear case demonstrates a different non-symmetry resulting from the incorporation of penalties.
The fiber layouts across all cases appear smooth, which can be attributed to the coarser 
B-spline discretization of the fiber orientation field relative to the state variable mesh.

 Comparison of the strain energy values in  Table \ref{cases_by_refinement} reveal slight variations 
 among the discretization orders, 
 with linear exhibiting the highest energy and cubic the lowest, 
 a consequence of the increased number of design variables in higher-order discretizations. 
These discretizations have slightly lower misalignment and curvature values compared to the linear case.

\begin{figure}[t]\centering
	\includegraphics[scale=1.0]{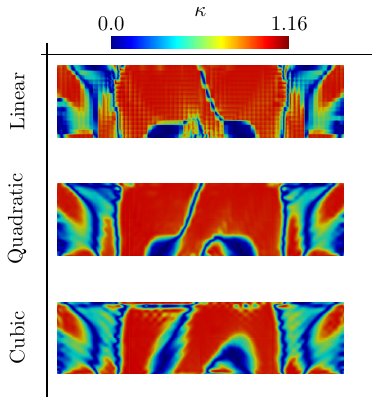}
	\caption{Curvature values for Case MLCP with varying polynomial order of the fiber orientation field discretization.}
	\label{fig_short_cantilver_beam_cubic_contour}
\end{figure}

 The higher-order discretization's smoothness becomes evident when evaluating the local curvature values which are 
 using derivatives of the fiber orientation field. 
 Fig.~\ref{fig_short_cantilver_beam_cubic_contour} shows local curvature values for Case MLCP across different 
 discretization orders, revealing discontinuities in the linear case within the contour plot, while quadratic and 
 cubic cases exhibit smoother curvature distributions.

This example illustrated the effectiveness of misalignment and curvature penalties in achieving smoother 
and more parallel fiber arrangements. It also revealed the influence of the polynomial order of fiber orientation field 
discretization on the optimized fiber layout. Although strain energy was largely unaffected, 
higher-order discretizations produce smoother first-order spatial derivatives which are subsequently used for calculating 
curvature and misalignment penalties.

\subsection{Cantilever beam under bending}\label{example2}

This example considers simultaneous optimization of the structural shape and fiber orientation in 2D. 
We consider quadratic polynomial order and two discretization sizes for the fiber orientation mesh, 
along with various penalization cases. 
This study extends the optimization problem to include the geometry level set field as a design parameter. Furthermore,
the influence of the initial fiber orientation on the final optimized layout is examined. 

\begin{figure}[t]\centering
	\includegraphics[scale=0.7]{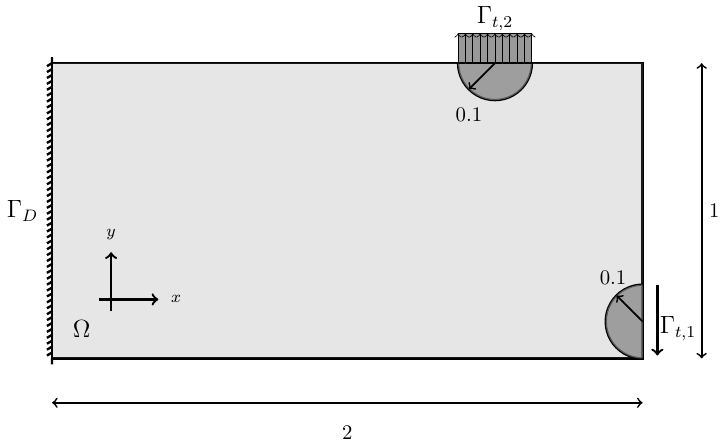}
	\caption{{Initial Configuration for concurrent fiber topology optimization, 2D.}}
	\label{fig_layout_exam2}
\end{figure}

\begin{figure*}[ht]\centering
	\includegraphics[scale=1.0]{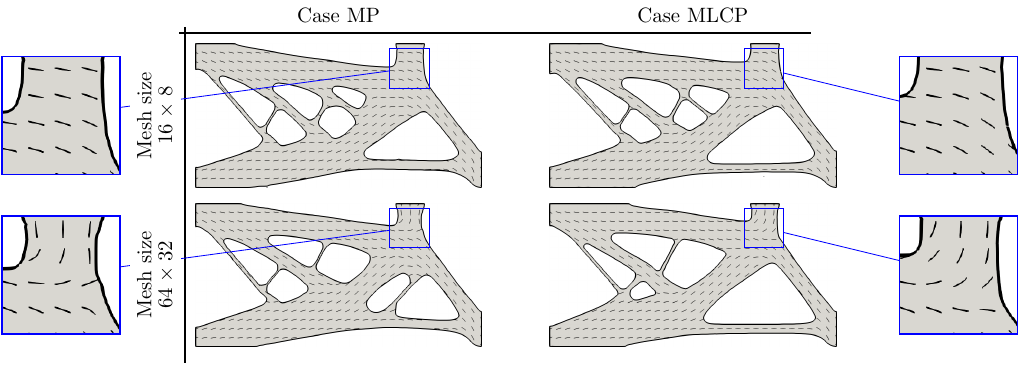}
	\caption{Streamline visualization of the concurrent geometry and fiber orientation optimized designs of a cantilever beam under bending for different mesh sizes and penalization cases.}
	\label{fig_concurrent_opt_design}
\end{figure*}

The objective is to minimize the compliance of the structure illustrated in Fig.~\ref{fig_layout_exam2}, 
which is subjected to two different loading cases under plane stress conditions. 
The cantilever beam is fixed on the left side, 
with non-design domains indicated by darker semicircles that transfer the loads. 
To ensure the fiber orientations are not solely dictated by principal stress directions, 
the problem is formulated with two loading scenarios. 
Loads are applied alternately at \( \Gamma_{t,1} \) and \( \Gamma_{t,2} \) with the values of \( P_{\Gamma_{t,1}} = 0.1 \) and \( P_{\Gamma_{t,2}} = 0.125 \). 
The optimization problem is formally expressed as follows:

\begin{equation} \label{Eq_example222}
	\begin{aligned}
	  \underset{\theta_{xy}, LS }{min} \Pi  = &w_{f} \left< \frac{\strainEnergy}{\strainEnergy^0} \right>  + \frac{w_p}{\mathcal{P}^0_p} \mathcal{P}_p + \frac{w_g}{\mathcal{P}^0_g}  \mathcal{P}_{g} \\ 
	  & \frac{w_{\text{par}}}{\mathcal{P}_{\text{par}}^{0} } \mathcal{P}_{\text{par}}  + \frac{w_{\text{Lcur}}}{\mathcal{P}_{\text{Lcur}}^{0} }  \mathcal{P}_{\text{Lcur}} \\ 
	  \mbox{s.t} & \quad V \leq \frac{1}{2}V_0.
	\end{aligned}
\end{equation}

Equation (\ref{Eq_example222}) components represent the averaged strain energy, alongside penalties for perimeter and 
regularization of the level set field introduced in Section \ref{RGLS}. The misalignment and local curvature penalties, as defined in Equation (\ref{eq:parallel_penalty2d3d}) 
and Equation (\ref{eq:curvature_penalty}), are also included. The operator 
$\left< \cdot \right>$ denotes the average of the strain energy over both load cases. 
The terms in the objective function are normalized by reference values denoted by superscript $0$.
The optimization is constrained by limiting the structural volume to no more than 50\% of the volume computational domain. 
The value of $\kappa_{\text{max}}$ is set to 10.0.

We consider two penalization cases: 
\begin{itemize}
	\item Case MP with only the parallel misalignment penalty ($w_{\text{par}} = 0.05$, $w_{\text{Lcur}} = 0.0$).
    \item Case MLCP with both misalignment and local curvature penalties ($w_{\text{par}} = 0.05$, $w_{\text{Lcur}} \ne 0.0$). 
\end{itemize}

State variables are discretized on a $256 \times 128$ mesh, while the level set field is discretized on a coarser $64 \times 32$ mesh. 
For the fiber orientation field, quadratic B-spline discretizations are utilized, with two and four times coarser meshes. 
This results in meshes with dimensions of $64 \times 32$ and $16 \times 8$ elements for the fiber orientation, respectively.

The initial level set field is created by seeding an array of holes arranged in a $7 \times 3$ grid, each with a radius of 0.12. 
The fiber orientation field is initialized to a constant value of zero, corresponding to horizontal fiber paths.

Figure \ref{fig_concurrent_opt_design} shows the optimized fiber layouts for various penalization cases and mesh sizes 
for the quadratic fiber orientation field. Table \ref{cantilverbeam_bending} shows the corresponding 
strain energy and penalty values, normalized against Case MP with a mesh size of $16 \times 8$, 
and the maximum curvature for each case.

\begin{table*}[ht]
    \centering
    \begin{tabular}{|c||c|c||c|c||}
    \hline
     & \multicolumn{2}{c||}{Case MP} & \multicolumn{2}{c||}{Case MLCP} \\
    \hline
     & $16 \times 8$ Mesh & $64 \times 32$ Mesh & $16 \times 8$ Mesh & $64 \times 32$ Mesh \\
    \hline
     $\strainEnergy$ & 1.000 & 0.890 & 1.091 & 0.907 \\
    \hline
    $max(\kappa)$ & 13.76 & 43.81 & 9.57 & 10.2 \\
    \hline
    $\mathcal{P}_{\text{par}}$ & 1.000 & 5.010 & 1.033 & 3.769 \\
    \hline
    $\mathcal{P}_{\text{Lcur}}$ & $1.000$ & $2.024\cdot 10^{6}$ & $0.000$ & $7.765\cdot 10^{-6}$ \\
    \hline
    \end{tabular}
	\caption{Objective components, normalized against Case MP with $16 \times 8$ Mesh size, and maximum curvature for different penalization cases and mesh sizes.}
	\label{cantilverbeam_bending}
\end{table*}

All cases in Fig.~\ref{fig_concurrent_opt_design} exhibit parallel fiber layouts, as anticipated due to the addition of the parallel penalty. 
The introduction of the curvature penalty results in smoother fiber paths. 
While this smoothness is not immediately evident for the configuration using a \(16 \times 8\) mesh, it becomes more pronounced when using a \(64 \times 32\) mesh. 
Additionally, the maximum curvature values observed in Table \ref{cantilverbeam_bending} 
for Case MLCP closely aligns with the maximum allowable curvature and is lower than the curvatures observed in Case MP. The 
high value of the curvature penalty in Case MP with the finer mesh size, \(64 \times 32\), is due to the kinks in the fiber path which results in 
a high curvature value. Strain energy values in Table \ref{cantilverbeam_bending} show that Case MLCP 
has higher strain energy than Case MP across different 
mesh sizes. This increase is due to the addition of another penalty term, which decreases the strain energy's contribution to the objective.

Investigating the influence of the discretization of the fiber orientation fields shows that a finer mesh leads to a smaller strain energy value. 
However, the coarser \(16 \times 8\) B-spline mesh leads to a smoother fiber orientation field and, consequently, 
 a more uniform fiber layout. The values for the parallel misalignment and curvature penalties are notably lower for the 
 coarser mesh compared to the finer one, indicating enhanced smoothness.

\begin{figure*}[!h]\centering
	\includegraphics[scale=1.0]{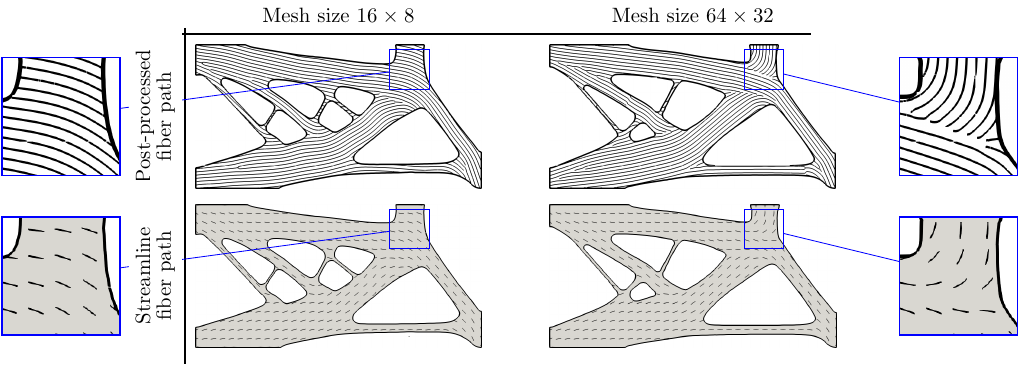}
	\caption{Post-processed and streamline visualization of the fiber paths for Case MLCP with quadratic B-spline discretization
	and two different mesh sizes.}
	\label{fig_print_ex2}
\end{figure*}

Figure \ref{fig_print_ex2} illustrates the post-processed fiber paths for the optimized designs for Case MLCP with two different mesh sizes, 
generated according to the process detailed in \cite{boddeti2020optimal} along with the streamline visualization. 
This figure illustrates that applying the parallel penalty reduces fiber intersections, while the curvature 
penalty minimizes fiber bending. Consequently, these penalties align the optimization results more closely with continuous, 
post-processed fiber paths.

\begin{figure*}[!h]\centering
	\includegraphics[scale=1.0]{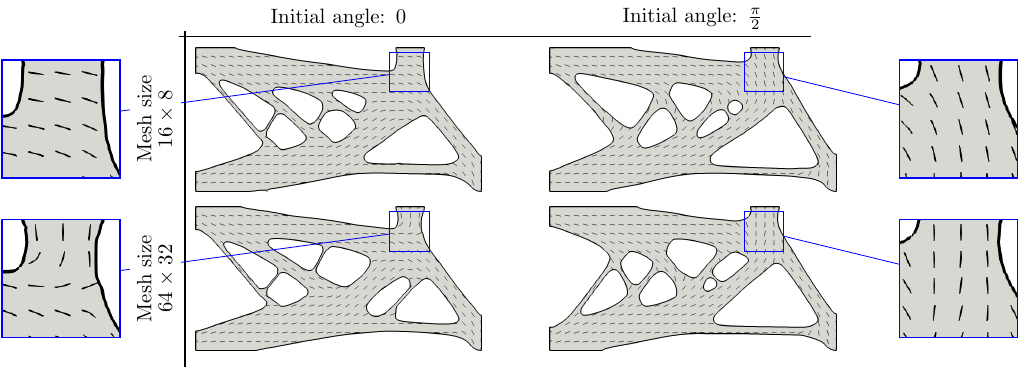}
	\caption{Streamline visualization of the concurrent geometry and fiber orientation optimized designs.}
	\label{fig_initail_dependency}
\end{figure*}

To assess the influence of the initial design of the fiber orientation field on the final optimized layout, 
the initial fiber layout is changed to vertical fiber paths, keeping the initial level set design unchanged. 
For Case MP and using a quadratic B-spline mesh, the optimization results are then compared to those obtained from designs initiated with horizontal fiber paths
in Fig.~\ref{fig_concurrent_opt_design}. Fig.~\ref{fig_initail_dependency} shows the optimized fiber layouts initialized with different angles. 
Table \ref{Initial dependency table} shows the corresponding strain energy 
and penalty values, normalized against the top left case in Fig.~\ref{fig_concurrent_opt_design}, as well as the absolute value of the maximum curvature for each case.

Designs initiated with horizontal fibers tend to retain more horizontal orientation in the optimized layout, 
and likewise, those starting with vertical fibers show a predominance of vertical fibers. 
When considering strain energy, it is observed that the values are slightly lower for designs 
initialized with vertical fibers than those starting with 
horizontal fibers. However, this difference is below 2\%. 
Thus, while having a significant impact on the fiber layout in some regions, the initial fiber orientation does not significantly 
affect the strain energy in this example.

\begin{table*}[ht]
    \centering
    \begin{tabular}{|c||c|c||c|c||}
    \hline
     & \multicolumn{2}{c||}{Initail angle: 0} & \multicolumn{2}{c||}{Initail angle : $\frac{\pi}{2}$} \\
    \hline
     & $16 \times 8$ Mesh & $64 \times 32$ Mesh & $16 \times 8$ Mesh & $64 \times 32$ Mesh \\
    \hline
     $\strainEnergy$ & 1.000 & 0.890 & 0.985 & 0.887 \\
    \hline
    $max(\kappa)$ & 13.76 & 43.81 & 16.3 & 61.2 \\
    \hline
    $\mathcal{P}_{\text{par}}$ & 1.000 & 5.010 & 1.757 & 5.183 \\
    \hline
    $\mathcal{P}_{\text{Lcur}}$ & $1.000\cdot 10^{0}$ & $2.024\cdot 10^{6}$ & $1.333\cdot 10^{4}$ & $3.184\cdot 10^{6}$ \\
    \hline
    \end{tabular}
	\caption{Objective components, normalized against the reference case, and maximum curvature for different initial fiber orientation angles.}
    \label{Initial dependency table}
\end{table*}

This example demonstrates the effectiveness of the misalignment and curvature penalties in generating more parallel and smoother fiber 
layouts for concurrent topology and fiber orientation optimization. 
The investigations into the parameterization of the fiber orientation field suggest the refinement level of the B-spline mesh plays an important role in the final optimized layout.
The results also highlighted the impact of the initial choice of fiber orientation on the final optimized layout.

\subsection{Support structure for a plate under uniform pressure}\label{example3}

This example minimizes the compliance of a 3D support structure for a plate,
 where both the structural shape and fiber orientations are considered design parameters. 
 The example examines the effects of misalignment and curvature penalties in cases where either a single fiber orientation, \(\theta_{xy}\), or both fiber 
 orientations, \(\theta_{xy}\) and \(\theta_z\), are treated as design variable fields.

\begin{figure}[b]\centering
	\includegraphics[scale=1.0]{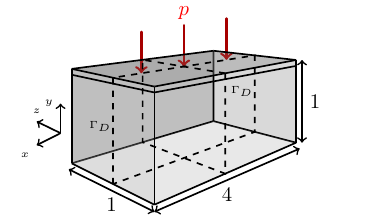}
	\caption{Support structure for a plate under uniform pressure.}
	\label{fig_3d_bridge}
\end{figure}

 The layout of the problem is illustrated in Fig.~\ref{fig_3d_bridge}, 
 showing the plate's top surface subjected to a uniform load of $p=0.01$. 
Dirichlet boundary conditions are applied to the left and right sides, 
 marked as $\Gamma_D$ in the figure.
The symmetry planes are shown by dashed lines. 
The plate, acting as the load-bearing component, is identified as the non-design domain ($0.95 \le z \le 1.0$),
 while the support structure is the design domain.

Utilizing the symmetry in the 
$x-y$ and $x-z$ planes, the simulation is conducted on a quarter of the structure. To ensure the symmetry of the fiber
orientation field, a penalty term is applied in the symmetry planes such that
on the $x-y$ symmetry plane $\theta_z = 0$ and on the $x-z$ symmetry plane $\theta_{xy} = 0$.
The quadratic mesh, used for discretizing the level set and fiber orientations, 
has dimensions of \(32 \times 16 \times 16\), making it twice as coarse as the linear mesh for state variables, which is \(128 \times 64 \times 64\).
 The optimization problem is formulated as follows:

\begin{equation} \label{Eq_examplebridge}
	\begin{aligned}
	  \underset{\theta_{xy}, (\theta_{z}), LS }{min} \Pi  = &w_{f} \frac{\strainEnergy}{\strainEnergy^0} + \frac{w_p}{\mathcal{P}_{p}^{0} } \mathcal{P}_p + \frac{w_g}{\mathcal{P}_{g}^{0} }  \mathcal{P}_{g} + \\ 
	  & \frac{w_{\text{par}}}{\mathcal{P}_{\text{par}}^{0} } \mathcal{P}_{\text{par}}+ \frac{w_{\text{Lcur}}}{\mathcal{P}_{\text{Lcur}}^{0} } \mathcal{P}_{\text{Lcur}} + \frac{w_{sym}}{\mathcal{P}_{sym}^{0} } \mathcal{P}_{sym}  \\ 
	  & \mbox{s.t}  \quad V \leq \frac{1}{4}V_0.
	\end{aligned}
\end{equation}

In Equation (\ref{Eq_examplebridge}), the terms represent the system's strain energy, perimeter and regularization penalties for the level set field, 
and the misalignment and local curvature penalty for fiber orientation field(s). The last penalty term enforces a 
symmetric fiber layout. Each term is normalized by reference values, 
indicated by the superscript $0$. 
The optimization is constrained to ensure the volume of the optimized design does not surpass 25\% of the computational domain's total volume. 
The maximum curvature is set to $5.0$. The initial level set field is an array of holes 
arranged in a $7 \times 3 \times 3 \times 3$ grid, each with a radius of 0.12, while the initial fiber
orientation field(s) is set to a constant value of zero, corresponding to horizontal fiber paths parallel to 
x-axis.

We define two cases to evaluate the effects of misalignment and local curvature penalty with different fiber orientation 
configurations as design variables:
\begin{itemize}
	\item Case 1: Considers solely $\theta_{xy}$ as the design variable field, while $\theta_z$ remains constant. A misalignment penalty with 
	a weight of 0.05 is applied. This setup mirrors the sequential, layer-wise printing process used in some FRC manufacturing techniques.
	\item Case 2: Incorporates both fiber orientations, $\theta_{xy}$ and $\theta_z$, as design variable fields, with a misalignment penalty weight of $0.05$.
\end{itemize}

\begin{figure*}[!h]\centering
	\includegraphics[scale=1.0]{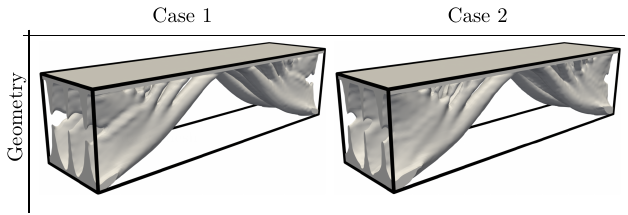}
	\caption{Optimized topology of the plate support structure for predefined cases.}
	\label{fig_3d_bridge_fiber_result}
\end{figure*}

\begin{figure*}[!h]\centering
	\includegraphics[scale=1.0]{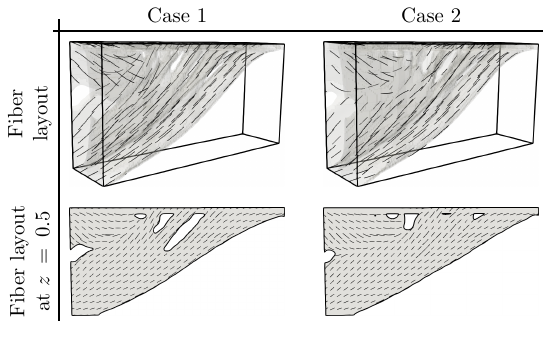}
	\caption{Streamline visualization of the optimized fiber layout and topology of the plate support structure for different cases.}
	\label{fig_3d_bridge_fiber_result2}
\end{figure*}

Figure~\ref{fig_3d_bridge_fiber_result} presents the optimized geometries for the two cases. 
The geometries display notable similarities, with slight variations primarily in the areas where the support structure connects to the plate. 
Figure~\ref{fig_3d_bridge_fiber_result2} illustrates the optimized topology and fiber layout for the specified cases, 
along with a cross-sectional slice at $z=0$. The fiber layout is similar for both cases, showing only slight variations in geometry within the slices.

Table \ref{EX3_fibers} presents the objective components, normalized against Case 1 and maximum curvature.
In both cases, the maximum curvature is close to the allowable limit of 5.0. 
In Case 2, controlling curvature becomes more challenging due to the potential alignment of fiber paths outside the \(x-y\) plane. 
This is reflected in the higher curvature penalty value and maximum curvature observed in this case.
It is noteworthy that the absolute curvature penalty is below \(10^{-9}\) for both cases ensuring the local curvature constraint is tightly satisfied.
Notably, Case 2 has a lower strain energy than Case 1, which is attributed to its expanded design space involving two design parameter fields
(\(\theta_{xy}\) and \(\theta_z\)), compared to Case 1's singular design parameter field (\(\theta_{xy}\)).

\begin{table}[!h]
	\centering
	\begin{tabular}{|c|c|c|c|c|}
	\hline
	Case & $\strainEnergy$ & max($\kappa$) & $\mathcal{P}_{par}$ & $\mathcal{P}_{Lcur}$ \\
	\hline
	1 & 1.000 & 5.15 & 1.000 & $1.000$ \\ \hline
	2 & 0.986 & 5.28 & 1.048 & $1.878\cdot 10^{6}$  \\ \hline
	\end{tabular}
	\caption{Objective components and maximum curvature for cases.}
	\label{EX3_fibers}
\end{table}

This example extended the misalignment and curvature penalties to 3D concurrent topology and fiber orientation optimization
and demonstrated their effectiveness in generating smoother and more parallel fiber layouts. Similar to the 2D configurations, 
the shape of the design domain was simple, i.e. it was rectangular.

\subsection{Cylinder with variable cross-section under torsion}\label{example4}

\begin{figure}[b]\centering
	\includegraphics[scale=3.0]{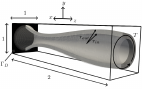}
	\caption{Cylindrical design domain with variable cross-section.}
	\label{fig_3d_vase}
\end{figure}

This example examines simultaneous fiber orientation and topology optimization for a complex geometric configuration,
featuring a design domain with a spatially variable circular 
cross-section and circumferential fiber orientation field.
The objective is to minimize strain energy, 
with both the level set and fiber orientation fields serving as design variables.

Figure~\ref{fig_3d_vase} illustrates the problem setup, where the axis-symmetric  design domain is  defined by an inner and an outer radius as follows:
\( r_{\mathrm{in}} = 0.15 + 0.2 z - 0.4 z^2 + 0.16 z^3 \), and \( r_{\mathrm{out}} = r_{\mathrm{in}} + 0.05 \). 
A non-design domain is placed at the cylinder's right end (\( 1.9 < z  < 2\)) to facilitate load transfer. 
A torque of $T = 0.05$ is uniformly applied over the cross-section at \(z=2\), while the left end of the cylinder is clamped.
This design domain is embedded into a rectangular computational domain. The geometry of the design domain is defined via the following level set fields:

\begin{equation}
	\begin{split}
		\phi_{1}(x,y,z) &= r_{\text{out}} - \sqrt{x^2 + y^2}, \\ 
		\phi_{2}(x,y,z) &= r_{\text{in}} - \sqrt{x^2 + y^2},
	\end{split}
\end{equation}

\noindent where \(r_{\text{out}}\) and \(r_{\text{out}}\) represent the outer and inner radii, respectively, and 
$(x,y,z)$ specify the spatial coordinates within the computational domain. 

A B-spline discretized design level set field is used to describe the shape of the solid within the design domain. 
Given that the fibers are arranged circumferentially, the angle $\theta_{xy}$ is obtained from the spatial coordinates using $\theta_{xy} = \atan2(y,x)$, 
hence it is not a design parameter. Conversely, $\theta_{z}$ is treated as a design parameter. 
A quadratic B-spline mesh of size  \(16 \times 16 \times 64\) is used to discretize the design level set field 
and fiber orientation field while a linear B-spline mesh of size \(64 \times 64 \times 256\) is used to discretize the state variable field.

The optimization problem is formally defined as:

\begin{equation} \label{Eq_example333}
	\begin{aligned}
	  \underset{\theta_{z}, LS }{min} \Pi  = &w_{f} \frac{\strainEnergy}{\strainEnergy^0}  + \frac{w_p}{\mathcal{P}^0_p} \mathcal{P}_p + \frac{w_g}{\mathcal{P}^0_g}  \mathcal{P}_{g} \\ 
	  & \frac{w_{\text{par}}}{\mathcal{P}_{\text{par}}^{0} } \mathcal{P}_{\text{par}} + \frac{w_{\text{Lcur}}}{\mathcal{P}_{\text{Lcur}}^{0} } \mathcal{P}_{\text{Lcur}}  \\ 
	  & \mbox{s.t}  \quad V \leq \frac{1}{2}V_0,
	\end{aligned}
\end{equation}

\noindent where the terms represent the system's strain energy, perimeter and regularization penalties for the level set, and misalignment and local curvature
penalties for the fiber orientation field. All terms in the objective function are normalized by the reference values indicated by the superscript $0$. 
A volume constraint is imposed, restricting the structural volume to 50\% of the design domain.
The weight for the parallel misalignment penalty is set to \( w_{\text{par}} = 0.05 \), and \(\kappa_{\text{max}}\) is assigned a value of 30.0. 
This value is higher than the values in the previous examples because the circumferential initialization of the fiber path means 
that the maximum curvature is \(max(\kappa) = \frac{1}{r_{\text{min}}} = 18.12\), with \(r_{\text{min}}\) 
being the minimum radius across the cross-section in the initial design. 

\begin{figure}[b]\centering
	\includegraphics[scale=3.0]{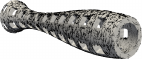}
	\caption{Initial design for the cylinder with variable cross-section.}
	\label{fig_3d_vase_init}
\end{figure}

Figure~\ref{fig_3d_vase_init} illustrates the initial design for both the geometry and fiber orientation fields.
The design level set field is initialized with a pattern of radially symmetric holes,
and the initial $\theta_{z}$ values are set to zero, resulting in
concentric fiber paths.

\begin{figure*}[!h]\centering
	\includegraphics[scale=3.0]{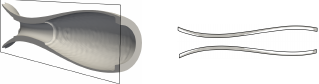}
	\caption{Optimized design geometry for the cylinder with variable cross-section.}
	\label{fig_3d_vase_final}
\end{figure*}

Figure~\ref{fig_3d_vase_final} shows the optimized design cropped with the plane $z=0.5$, 
and the corresponding cross-section. Since the optimization aims to minimize compliance, it effectively seeks to 
maximize the torsional rigidity of the structure by creating a cross-section that has a nearly uniform torsional stiffness along the length of the cylinder.
 For a thin circular cross-section, 
torsional rigidity is described by $J=\frac{2}{3} \pi r^3t$, with $r$ the radius and $t$ thickness. 
The figure demonstrates that the variation in thickness along the z-axis is inversely proportional to the radius,
consistent with the torsional rigidity equation in a thin circular cross-section.
Consequently, in areas with a larger radius, the structure is thinner, whereas in the central region with a smaller radius,
 the thickness increases, occupying more of the design domain.

\begin{figure*}[!h]\centering
	\includegraphics[scale=1.0]{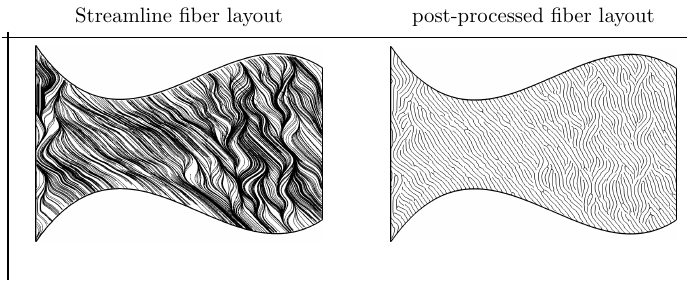}
	\caption{Streamline and post-processed fiber layout for the optimized design of the cylinder with variable cross-section.}
	\label{fig_3d_vase_fibers_layout}
\end{figure*}

Figure~\ref{fig_3d_vase_fibers_layout} shows the optimized fiber layout visualized with streamlines alongside the post-processed 
continuous post-processed fiber layout. 
The fiber layouts are depicted by projecting the design's outer boundary surface onto a plane. 
Note that since a surface with non-zero Gaussian curvature is mapped into a plane, distortion increases 
progressively away from the central line of the plane.

In a state of uniform torsion, the outer boundary elements of the structure are subjected to pure shear, 
resulting in principal stress directions to align at $\frac{\pi}{4}$ and $\frac{3\pi}{4}$. 
Since this problem is shear-dominated,  
fiber angles alternate between orientations near $\frac{\pi}{4}$ and $\frac{3\pi}{4}$ with transition 
regions between them.
The figure also indicates that the penalty terms assist in aligning and smoothing the optimized fiber paths,
 making them more closely resemble the post-processed continuous fiber paths.

The absolute values of the objective components and maximum curvature for the optimized design are presented in Table \ref{EX4_ref_plate}. 
As it can be seen, the maximum curvature is close to the allowable limit of 30.0, and the local curvature penalty is below $10^{-5}$, 
indicating that the local curvature constraints imposed by the local curvature penalty are tightly satisfied.

\begin{table}[ht]
    \centering
    \begin{tabular}{|c|c|}
    \hline
    Parameter & Value \\\hline
    $\strainEnergy$ & $5.18\cdot 10^{-6}$ \\ \hline
    max($\kappa$) & $30.3$ \\ \hline
    $\mathcal{P}_{par}$ & $9.646\cdot 10^{-2}$ \\ \hline
    $\mathcal{P}_{Lcur}$ & $3.542\cdot 10^{-6}$ \\ \hline
    \end{tabular}
    \caption{Absolute values of objective components and maximum curvature for the optimized design of the cylinder with variable cross-section.}
	\label{EX4_ref_plate}
\end{table}